\documentclass[12pt]{article}

\usepackage{amsmath, amssymb}

\numberwithin{equation}{section}
\newcommand{\nop}{\mbox{$\circ\atop\circ$}}
\newcommand{\1}{\mathbf {1}}
\newcommand{\Z}{{\mathbb Z}}
\newcommand{\Q}{{\mathbb Q}}
\newcommand{\R}{{\mathbb R}}
\newcommand{\C}{{\mathbb C}}
\newcommand{\g}{{\mathfrak g}}
\newcommand{\h}{{\mathfrak h}}
\newcommand{\CH}{{\mathcal H}}
\newcommand{\I}{{\mathcal I}}
\newcommand{\CP}{{\mathcal P}}
\newcommand{\W}{{\mathcal W}}

\newcommand{\widet}{\widetilde}
\newcommand{\al}{\alpha}
\newcommand{\be}{\beta}

\newcommand{\gm}{\gamma}
\newcommand{\om}{\omega}
\newcommand{\sg}{\sigma}
\newcommand{\ve}{\varepsilon}
\newcommand{\prf}{\noindent {\bfseries Proof} \  }
\newcommand{\qed}{\mbox{ $\square$}}
\newcommand{\la}{\langle}

\newcommand{\ot}{\otimes}
\newcommand{\ra}{\rangle}
\newcommand{\bv}{\mathbf{v}}

\DeclareMathOperator{\Aut}{Aut}
\DeclareMathOperator{\End}{End}
\DeclareMathOperator{\Res}{Res}
\DeclareMathOperator{\Vir}{Vir}
\DeclareMathOperator{\ch}{ch}
\DeclareMathOperator{\spn}{span}
\DeclareMathOperator{\tr}{tr}
\DeclareMathOperator{\wt}{wt}

\newtheorem{thm}{Theorem}[section]
\newtheorem{prop}[thm]{Proposition}
\newtheorem{lem}[thm]{Lemma}
\newtheorem{cor}[thm]{Corollary}
\newtheorem{rmk}[thm]{Remark}

\begin{document}

\begin{center}
{\Large {\bfseries $\Z_3$ symmetry and $W_3$ algebra in lattice
vertex operator algebras}}

\vspace{10mm}

Chongying Dong\footnote{Partially
supported by NSF grant DMS-9987656 and a research grant from the
Committee on Research, UC Santa Cruz.}\\
Department of Mathematics, University of California\\
Santa Cruz, CA 95064\\

\vspace{2mm} Ching Hung Lam\footnote{Partially supported by NSC
grant 91-2115-M-006-014 of Taiwan, R.O.C.}\\
Department of Mathematics, National Cheng Kung University\\
Tainan, Taiwan 701\\

\vspace{2mm} Kenichiro Tanabe\footnote{Partially supported by JSPS
Grant-in-Aid for Scientific Research No. 14740061.}\\
Institute of Mathematics, University of Tsukuba\\
Tsukuba 305-8571, Japan\\

\vspace{2mm}

Hiromichi Yamada\footnote{Partially supported by JSPS
Grant-in-Aid for Scientific Research No. 13640012.}\\
Department of Mathematics, Hitotsubashi University\\
Kunitachi, Tokyo 186-8601, Japan\\

\vspace{2mm}
Kazuhiro Yokoyama\\
Graduate School of Mathematics, Kyushu University\\
Fukuoka 812-8581, Japan
\end{center}

\vspace{10mm}
\begin{abstract}
The $W_3$ algebra of central charge $6/5$ is realized
as a subalgebra of the vertex operator algebra
$V_{\sqrt{2}A_2}$ associated with a lattice of
type $\sqrt{2}A_2$ by using both coset construction and
orbifold theory. It is proved that $W_3$ is rational.
Its irreducible modules are
classified and constructed explicitly. The
characters of those irreducible modules are also computed.
\end{abstract}

\section{Introduction}
The vertex operator algebras associated with positive definite
even lattices affords a large family of known examples of vertex
operator algebras. An isometry of the lattice induces an automorphism
of the lattice vertex operator algebra.
The subalgebra of fixed points is the so called orbifold vertex
operator algebra. In this paper we deal with the case where the lattice
$L=\sqrt{2}A_2$ is
$\sqrt{2}$ times an ordinary root lattice of type $A_2$ and the
isometry $\tau$ is an element of the Weyl group of order $3$.
We use this algebra to study the $W_3$ algebra of central charge $6/5.$
In fact, by using both coset construction and orbifold theory we construct
the $W_3$ algebra of central charge
$6/5$ inside $V_L$ and classify its irreducible modules. We also prove
that the $W_3$ algebra
is rational and compute the characters of the irreducibles.

The vertex operator algebra $V_L$
associated with $L=\sqrt{2}A_2$ contains three mutually orthogonal
conformal vectors $\om^1$, $\om^2$, $\om^3$ with central charge
$c=1/2$, $7/10$, or $4/5$ respectively \cite{DLMN}. The subalgebra
$\Vir(\om^i)$ generated by $\om^i$ is the Virasoro vertex operator
algebra $L(c,0)$, which is the  irreducible unitary highest
weight module for the Virasoro algebra with central charge $c$ and
highest weight $0$. The structure of $V_L$ as a module for
$\Vir(\om^1) \otimes \Vir(\om^2) \otimes \Vir(\om^3)$ was
discussed in \cite{KMY}. Among other things it was shown that
$V_L$ contains a subalgebra of the form $L(4/5,0)\otimes
L(4/5,3)$. Such a vertex operator algebra is called a $3$-state
Potts model. This subalgebra is contained in the subalgebra
$(V_L)^\tau$ of fixed points of $\tau$. There is another
subalgebra $M$ in $V_L$, which is of the form
\begin{equation*}
L(\frac{1}{2},0) \otimes L(\frac{7}{10},0) \oplus
L(\frac{1}{2},\frac{1}{2}) \otimes L(\frac{7}{10},\frac{3}{2})
\end{equation*}
and is invariant under $\tau$. The representation theory of $M$
was studied in \cite{KLY1, LY}.

We are interested in the subalgebra $M^\tau$ of fixed points of
$\tau$ in $M$. Its Virasoro element is $\om=\om^1+\om^2$. The
central charge of $\omega$ is $1/2 + 7/10 = 6/5$. We find an
element $J$ of weight $3$ in $M^\tau$ such that the component
operators $L(n)=\om_{n+1}$ and $J(n)=J_{n+2}$ satisfy the same
commutation relations as in \cite[(2.1), (2.2)]{BMP} for $W_3.$
Thus the vertex operator subalgebra $\W$ generated by $\om$ and $J$ is a
$W_3$ algebra with central charge $6/5$.

We construct $20$ irreducible $M^\tau$-modules. $8$ of them are
inside irreducible untwisted $M$-modules, while $6$ of them are
inside irreducible $\tau$-twisted $M$-modules and the remaining
$6$ are inside irreducible $\tau^2$-twisted $M$-modules. There are
exactly two inequivalent irreducible $\tau^i$-twisted $M$-modules
$M_T(\tau^i)$ and $W_T(\tau^i)$, $i=1,2$. We investigate the
irreducible $\tau^i$-twisted $V_L$-modules constructed in
\cite{DL} and obtain $M_T(\tau^i)$ and $W_T(\tau^i)$ inside them.

We classify the irreducible modules for $\W$ by
determining the Zhu algebra $A(\W)$ (cf. \cite{Z}).
The method used here is similar to that in \cite{Wa2}, where the
Zhu algebra of a $W_3$ algebra with central charge $-2$ is
studied. We can define a map of the polynomial algebra $\C [x,y]$
with two variables $x, y$ to $A(\W)$ by $x \mapsto [\omega]$ and
$y \mapsto [J]$, which is a surjective algebra homomorphism. Thus
it is sufficient to determine its kernel $\I$. The key point is
the existence of a singular vector $\bv$ for the $W_3$ algebra
$\W$ of weight $12$. An invariant positive definite hermitian form
on $V_L$ implies that $\bv$ is in fact $0$. Thus $[\bv] = 0$.
Moreover, $[J(-1)\bv]=[J(-2)\bv]=[J(-1)^2\bv]=0$. Hence the
corresponding polynomials in $\C [x,y]$ must be contained in the
ideal $\I$. It turns out that $\I$ is generated by those four
polynomials and the classification of irreducible $\W$-modules is
established by Zhu's theory (\cite{Z}). That is, there are exactly
$20$ inequivalent irreducible $\W$-modules. The calculation of
explicit form of the singular vector $\bv$ and the calculation of
the ideal $\I$ were done by a computer algebra system Risa/Asir.

By the classification of irreducible $\W$-modules and an invariant
positive definite hermitian form, we can show that $M^{\tau} =
\W$. The eigenvalues of the action of weight preserving operators
$L(0)=\om_1$ and $J(0)=J_2$ on the top levels of those $20$
irreducible $M^\tau$-modules coincide with the values $\Delta
\begin{pmatrix} n & m\\ n'& m'\end{pmatrix}$ and $w
\begin{pmatrix} n & m\\ n'& m'\end{pmatrix}$ of \cite[(1.2),
(5.6)]{FZ} with $p=5$. Hence our $M^{\tau}$ is an algebra denoted
by $[Z_3^{(5)}]$ in \cite{FZ}.

We prove that $\W$ is $C_2$-cofinite and rational by using the
singular vector $\bv$ of weight $12$ and the irreducible modules
for $\W.$ In the course of proof we obtain a result about a
general vertex operator algebra $V.$ It says that if $V$ is
$C_2$-cofinite, then $V$ is rational if and only if $A(V)$ is
semisimple and any simple $A(V)$-module generates an irreducible
$V$-module. This result itself is not very hard to prove. But it
will certainly be useful in the future study of relationship
between rationality and $C_2$-cofiniteness.

We also study the characters of those irreducible
$M^{\tau}$-modules. Using the modular invariance of trace
functions in orbifold theory (cf. \cite{DLM}),
we describe the characters of
the $20$ irreducible $M^{\tau}$-modules in terms of the characters
of irreducible unitary highest weight modules for the Virasoro
algebras.

The results in this paper have applications to the Monster simple group.
Recently, it
was shown in \cite{KLY2} that the $\Z_3$ symmetry of a $3$-state
Potts model in $(V_L)^\tau$ affords $3A$ elements of the Monster
simple group. Such a result has been suggested by \cite{M2}. It is
expected that the $\Z_3$ symmetry of $M^{\tau}$ affords $3B$
elements.

The organization of the paper is as follows. In Section 2 we
review some properties of $M$ for later use. In Section 3 we
define the vector $J$ and compute the commutation relations among
the component operators $L(n)=\om_{n+1}$ and $J(n)=J_{n+2}$. In
Section 4 we construct $20$ irreducible $M^\tau$-modules and
discuss their properties. In Section 5 we determine the Zhu
algebra of the vertex operator subalgebra $\W$ generated by
$\omega$ and $J$ and show that $M^{\tau} = \W$. Thus we conclude
that $M^{\tau}$ has exactly $20$ inequivalent irreducible modules.
Finally, in Section 6 we study the characters of those irreducible
$M^{\tau}$-modules.

The authors would like to thank Toshiyuki Abe and Kiyokazu
Nagatomo for helpful advice concerning $W$ algebras and Hiroshi
Yamauchi for comment on an invariant positive definite hermitian
form on $V_L$. They also would like to thank Masahiko Miyamoto for
valuable discussions.

\section{Subalgebra $M$ of $V_{\sqrt{2}A_2}$}

In this section we fix notation. For basic definitions
concerning lattice vertex operator algebras we refer to \cite{
DL, FLM}. We also recall certain properties of the vertex
operator algebra $V_{\sqrt{2}A_2}$ (cf. \cite{KMY}).

Let $\al_1, \al_2$ be the simple roots of type $A_2$ and set
$\al_0=-(\al_1+\al_2)$. Then $\la \al_i, \al_i\ra = 2$ and $\la
\al_i, \al_j\ra=-1$ if $i \ne j$. Set $\be_i=\sqrt{2}\al_i$ and
let $L=\Z\be_1+\Z\be_2$ be the lattice spanned by $\be_1$ and
$\be_2$. We usually denote $L$ by $\sqrt{2}A_2$.

We follow Sections 2 and 3 of \cite{DL} with $L=\sqrt{2}A_2$,
$p=3$, and $q=6$. In our case $\la\al,\be\ra \in 2\Z$ for all
$\al, \be \in L$, so that the alternating $\Z$-bilinear map $c_0 :
L \times L \to \Z/6\Z$ defined by \cite[(2.9)]{DL} is trivial.
Thus the central extension
\begin{equation}\label{ext1}
1 \longrightarrow \la\kappa_6\ra \longrightarrow \widehat{L}
\overset{-}{\longrightarrow} L \longrightarrow 1
\end{equation}
determined by the commutator condition
$aba^{-1}b^{-1}=\kappa_6^{c_0(\bar{a},\bar{b})}$ splits. Then for
each $\al \in L$, we can choose an element $e^\al$ of $\hat{L}$ so
that $e^\al e^\beta=e^{\al+\beta}$. The twisted group algebra
$\C\{L\}$ is isomorphic to the ordinary group algebra $\C[L]$.

We adopt the same notation as in \cite{KLY1} to denote cosets of
$L$ in the dual lattice $L^\perp = \{ \al \in \Q \ot_{\Z}
L\,|\,\la\al, L\ra \subset \Z\}$, namely,

\begin{equation*}
L^0=L,\quad  L^1=\frac{-\be_1+\be_2}{3}+L ,\quad
L^2=\frac{\be_1-\be_2}{3}+L,
\end{equation*}
\begin{equation*}
L_0=L,\quad L_a=\frac{\be_2}{2}+L,\quad
L_b=\frac{\be_0}{2}+L,\quad L_c=\frac{\be_1 }{2}+L,
\end{equation*}
and
\begin{equation*}
L^{(i,j)} = L_i + L^j
\end{equation*}
for $i=0,a,b,c$ and $j=0,1,2$, where $\{0,a,b,c\} \cong \Z_2
\times \Z_2$. Then, $L^{(i,j)},i \in \{0,a,b,c\},j\in \{0,1,2\}$
are all the cosets of $L$ in $ L^{\perp }$ and $L^\perp/L \cong
\Z_2 \times \Z_2 \times \Z_3$.

Our notation for the vertex operator algebra $(V_L,
Y(\,\cdot\,,z))$ associated with $L$ is standard \cite{FLM}. In
particular, ${\mathfrak h}=\C\otimes_{\Z} L$  is an abelian Lie
algebra, $\hat {\mathfrak h}={\mathfrak h}\otimes
\C[t,t^{-1}]\oplus \C c$ is the corresponding affine Lie algebra,
$M(1)=\C[\al(n)\,;\,\al\in {\mathfrak h}, n<0],$ where
$\al(n)=\al\otimes t^n,$ is the unique irreducible $\hat{\mathfrak
h}$-module such that $\alpha(n)1=0$ for all $\alpha\in {\mathfrak
h}$ and $n>0$, and $c=1$. As a vector space $V_L = M(1) \ot \C[L]$
and for each $v \in V_L$, a vertex operator $Y(v, z) = \sum_{n \in
\Z} v_n z^{-n-1} \in \End (V_L)[[z,z^{-1}]]$ is defined. The
coefficient $v_n$ of $z^{-n-1}$ is called a component operator.
The vector $\1 = 1 \otimes 1$ is called the vacuum vector.

By Dong \cite{D1}, there are exactly $12$ isomorphism classes of
irreducible $V_L$-modules, which are represented by
$V_{L^{(i,j)}}$, $i=0,a,b,c$ and $j=0,1,2$. We use the symbol
$e^{\alpha}, \alpha \in L^\perp$ to denote a basis of $\C\{
L^\perp\}$.

To describe certain weight $2$ elements in $V_L$, we introduce the
following notation.
\begin{equation*}
x(\al)=e^{\sqrt{2}\al} + e^{-\sqrt{2}\al},\qquad
y(\al)=e^{\sqrt{2}\al} - e^{-\sqrt{2}\al},\qquad
w(\al)=\frac{1}{2}\al(-1)^2 - x(\al)
\end{equation*}
for $\al \in\{\pm\al_0, \pm\al_1, \pm\al_2\}$. We have
\begin{equation}\label{w1w}
w(\al_i)_1 w(\al_j) =
\begin{cases} 8w(\al_i) & \text{if } i=j\\
w(\al_i) + w(\al_j) - w(\al_k) & \text{if } i\ne j,
\end{cases}
\end{equation}
where $k$ is such that $\{i,j,k\}=\{0,1,2\}$. Moreover,
$w(\al_i)_2 w(\al_j) =0$ and
\begin{equation}\label{w3w}
w(\al_i)_3 w(\al_j) =
\begin{cases} 4\1 & \text{if } i=j\\
\frac{1}{2}\1 & \text{if } i\ne j.
\end{cases}
\end{equation}

Let
\begin{gather*}
\om = \frac{1}{5}\big( w(\al_1) + w(\al_2) + w(\al_0)\big),\\
\widet{\om} = \frac{1}{6}\big( \al_1(-1)^2 + \al_2(-1)^2 +
\al_0(-1)^2\big),\\
\om^1 = \frac{1}{4}w(\al_1),\qquad \om^2 = \om - \om^1,\qquad
\om^3 = \widet{\om} - \om.
\end{gather*}
Then $\widet{\om}$ is the Virasoro element of $V_L$ and $\om^1,
\om^2,\om^3$ are mutually orthogonal conformal vectors of central
charge $1/2, 7/10, 4/5$ respectively (cf. \cite{DLMN}). The
subalgebra $\Vir(\om^i)$ generated by $\om^i$ is isomorphic to the
Virasoro vertex operator algebra of given central charge, and
$\om^1$, $\om^2$, and $\om^3$ generate
\begin{equation*}
\Vir(\om^1) \otimes \Vir(\om^2) \otimes \Vir(\om^3) \cong
L(\frac{1}{2},0) \otimes L(\frac{7}{10},0) \otimes
L(\frac{4}{5},0).
\end{equation*}

We study certain subalgebras, and also submodules for them
in $V_{L_i}, i=0,a,b,c$ and $V_{L^j}, j=0,1,2$. Set
\begin{align*}
M_{k}^{i}&=\{ v\in V_{L_i}\,|\,(\omega^3)_1 v=0\},\\
W_{k}^{i}&=\{ v\in V_{L_i}\,|\,(\omega^3)_1 v=
\frac{2}{5}v\} ,\qquad \text{ for } i=0,a,b,c,
\end{align*}
and
\begin{align*}
M_{t}^{j}&=\{ v\in V_{L^j}\,|\,(\omega^1)_1 v=
(\omega^2)_1 v=0\},\\
W_{t}^{j}&=\{ v\in V_{L^j}\,|\,(\omega^1)_1 v=0,\quad
(\omega^2)_1 v=\frac{3}{5}v\} ,\qquad
\text{ for } j=0,1,2.
\end{align*}

Then $M_k^0$ and $M_t^0$ are simple vertex operator algebras.
Furthermore, $\{M_k^i,\,W_k^i,\, i=0,a,b,c\}$ and $\{M_t^j,\,
W_t^j,\, j=0,1,2\}$ are the sets of all inequivalent irreducible
modules for $M_k^0$ and $M_t^0$, respectively
\cite{KLY1, KMY,LY}. We also have

\begin{align*}
& M_{k}^{0}\cong L( \frac{1}{2},0) \otimes L( \frac{7}{10} ,0)
\oplus L( \frac{1}{2},\frac{1}{2}) \otimes L
\frac{7}{10},\frac{3}{2}),\\
& W_{k}^{0}\cong L( \frac{1}{2},0)
\otimes L( \frac{7}{10},\frac{3}{5}) \oplus L( \frac{1}{2},
\frac{1}{2}) \otimes L(\frac{7}{10},\frac{1}{10}),\\
& M_{k}^{a}\cong M_{k}^{b}\cong L( \frac{1}{2},\frac{1}{16})
\otimes L(\frac{7}{10},\frac{7}{16}),\\
& W_{k}^{a}\cong W_{k}^{b}\cong L( \frac{1}{2},\frac{1}{16})
\otimes L( \frac{7}{10},\frac{3}{80}),\\
& M_{k}^{c}\cong L( \frac{1}{2},\frac{1}{2}) \otimes
L(\frac{7}{10},0) \oplus L( \frac{1}{2},0) \otimes
L(\frac{7}{10},\frac{3}{2}),\\
& W_{k}^{c}\cong L( \frac{1}{2},\frac{1}{2}) \otimes L(
\frac{7}{10},\frac{3}{5}) \oplus L( \frac{1}{2},0) \otimes
L(\frac{7}{10},\frac{1}{10})
\end{align*}
as $L(1/2,0)\otimes L(7/10,0)$-modules and
\begin{align*}
M_{t}^{0} & \cong L( \frac{4}{5},0) \oplus L( \frac{4}{5},3), &
W_{t}^{0} & \cong L( \frac{4}{5},\frac{2}{5})
\oplus L( \frac{4}{5},\frac{7}{5}),\\
M_{t}^{1} & \cong M_{t}^{2}\cong L( \frac{4}{5},\frac{2}{3}), &
W_{t}^{1} & \cong W_{t}^{2}\cong L( \frac{4}{5},\frac{1}{15})
\end{align*}
as $L(4/5,0)$-modules.

Note also that
\begin{equation*}
V_{L^{(i,j)}}\cong \left(M_k^i\otimes M_t^j\right) \oplus
\left(W_k^i\otimes W_t^j\right)
\end{equation*}
as an $M_k^0\otimes M_t^0$-module.

We consider the following three isometries of
$(L,\la\cdot,\cdot\ra)$:
\begin{align*}
\tau &: \be_1 \to \be_2 \to \be_0 \to \be_1,\\
\sg &: \be_1 \to \be_2, \qquad \be_2 \to \be_1,\\
\theta &: \be_i \to -\be_i, \quad i=1,2.
\end{align*}

Note that $\tau$ is fixed-point-free and of order $3$. Note also
that $\sg \tau \sg = \tau^{-1}$. The isometries $\tau, \sg$, and
$\theta$ of $L$ can be extended to isometries of $L^\perp$. Then
they induce permutations on $L^\perp/L$. Since $\hat{L}$ is a
split extension, the isometry $\tau$ of $L$ lifts naturally to an
automorphism of $\hat{L}$. Then it induces an automorphism of
$V_L$:
\begin{equation*}
\al^1(-n_1)\cdots\al^k(-n_k)e^{\be} \longmapsto
(\tau\al^1)(-n_1)\cdots(\tau\al^k)(-n_k)e^{\tau\be}.
\end{equation*}

By abuse of notation, we  denote it by $\tau$ also. Moreover,
we can consider the action of $\tau$ on $V_{L^{(i,j)}}$ in a
similar way. We apply the same argument to $\sigma$ and $\theta$.

Set $M=M_k^0$. The vertex operator algebra $M$ plays an important
role in this paper. Recall that
\begin{equation*}
M \cong L(\frac{1}{2},0) \ot L(\frac{7}{10},0) \oplus
L(\frac{1}{2},\frac{1}{2}) \ot L(\frac{7}{10},\frac{3}{2})
\end{equation*}
as $\Vir(\om^1)\otimes \Vir(\om^2)$-modules. Note that $\om$ is
the Virasoro element of $M$ whose central charge is $6/5$. For $u
\in M$, we have $\om_1 u = hu$ for some $h \in \Z$ if and only if
$\widetilde{\om}_1 u=hu$. In such a case $h$ is called the weight
of $u$. Note also that $M$ is generated by $w(\al_1)$, $w(\al_2)$,
and $w(\al_0)$. In particular, $M$ is invariant under $\tau$,
$\sigma$, and $\theta$. In fact, $\theta$ acts on $M$ as the
identity.

We next show that the automorphism group $\Aut(M)$ of $M$ is
generated by $\sigma$ and $\tau$.

\begin{thm}
$(1)$ There are exactly three conformal vectors of central charge
$1/2$ in $M$, which are $\frac{1}{4}w(\al_i)$, $i=0,1,2$.

$(2)$ $\Aut(M)=\la\sg,\tau\ra$ is isomorphic to a symmetric group
of degree $3$.
\end{thm}

\prf We first consider conformal vectors in $M$. By \cite[Lemma
5.1]{M1}, a weight $2$ vector $v$ is a conformal vector of central
charge $1/2$ if and only if $v_1 v=2v$ and $v_3 v=\frac{1}{4}\1$.
Since $\{w(\al_0), w(\al_1), w(\al_2)\}$ is a basis of the weight
$2$ subspace of $M$, we may write $v = \sum_{i=0}^2 a_i w(\al_i)$
for some $a_i \in \C$. From (\ref{w1w}) and (\ref{w3w}) we see
that $v_1 v=2v$ and $v_3 v=\frac{1}{4}\1$ hold only if
$(a_0,a_1,a_2)=(1/4,0,0)$, $(0,1/4,0)$, or $(0,0,1/4)$. This
proves (1). Then any automorphism of $M$ induces a permutation on
$\{w(\al_0), w(\al_1), w(\al_2)\}$. If an automorphism induces the
identity permutation on the set, it must be the identity since $M$
is generated by $w(\al_1)$, $w(\al_2)$, and $w(\al_0)$. Now
\begin{equation*}
\tau : w(\al_1) \to w(\al_2) \to w(\al_0) \to w(\al_1),
\end{equation*}
and \begin{equation*} \sigma : w(\al_1) \to w(\al_2),\quad
w(\al_2) \to w(\al_1),\quad w(\al_0) \to w(\al_0).
\end{equation*}
Hence (2) holds. \quad \qed

\bigskip Let $v_h = w(\al_2) - w(\al_0)$. This vector is a highest
weight vector of highest weight $(1/2, 3/2)$ for $\Vir(\om^1)
\otimes \Vir(\om^2)$, that is, $(\om^1)_1 v_h = (1/2)v_h$,
$(\om^2)_1 v_h = (3/2)v_h$, and $(\om^1)_n v_h = (\om^2)_n v_h =
0$ for $n \ge 2$. Thus the $\Vir(\om^1) \otimes
\Vir(\om^2)$-submodule in $M$ generated by $v_h$ is isomorphic to
$L(1/2, 1/2)\otimes L(7/10,3/2)$. In particular, $M$ is generated
by $\om^1$, $\om^2$, and $v_h$.

We can choose another generator of $M$. Let
\begin{equation}\label{u1u2}
u^1 = w(\al_1)+\xi^2 w(\al_2)+\xi w(\al_0), \qquad
u^2=w(\al_1)+\xi w(\al_2)+\xi^2 w(\al_0),
\end{equation}
where $\xi = \exp(2\pi\sqrt{-1}/3)$ is a primitive cubic root of
unity. Then $\tau u^1=\xi u^1$, $\tau u^2=\xi^2 u^2$, and $\sigma
u^1=\xi^2 u^2$. We also have $(u^1)_1 u^1 = 4u^2$ and $((u^1)_1
u^1)_1 u^1=140\om$. Thus $u^1$, $(u^1)_1 u^1$, and $((u^1)_1
u^1)_1 u^1$ span the weight $2$ subspace of $M$. This implies that
$M$ is generated by a single vector $u^1$. A similar assertion
holds for $u^2$.

The subalgebra $M_t^0\cong L(4/5,0) \oplus L(4/5,3)$ is called a
$3$-state Potts model. It plays an important role in Subsection
4.2. The irreducible $M_t^0$-modules and their fusion rules are
determined in \cite{KMY, M2}. The Virasoro element of $M_t^0$ is
$\om^3$. Let
\begin{equation}\label{hwvt}
  \begin{split}
  v_t=&\frac{1}{9} (\al_1-\al_2)(-1)(\al_2-\al_0)(-1)(\al_0-\al_1)(-1)\\
  &-\frac{1}{2} (\al_1-\al_2)(-1)x(\al_0)
  -\frac{1}{2} (\al_2-\al_0)(-1)x(\al_1)
  -\frac{1}{2} (\al_0-\al_1)(-1)x(\al_2),
  \end{split}
\end{equation}
which is denoted by $q$ in \cite{KMY}. The vector $v_t$ is a
highest weight vector in $M_t^0$ of highest weight $3$ for
$\Vir(\om^3)$. Clearly, $\tau v_t = v_t$ and thus $\tau$ fixes
every element in $M_t^0$. Moreover, $\sigma v_t = -v_t$ and
$\theta v_t=-v_t$. Hence $\sigma$ and $\theta$ induce the same
automorphism of $M_t^0$, namely, $1$ on $\Vir(\om^3) \cong
L(4/5,0)$ and $-1$ on the $\Vir(\om^3)$-submodule generated by
$v_t$, which is isomorphic to $L(4/5,3)$. The automorphism group
$\Aut (M_t^0)$ is of order $2$ generated by $\theta$.

\section{Subalgebra $\W$ generated by $\om$ and $J$ in $M^\tau$}

For any $\tau$-invariant space $U$, set $U(\epsilon)=\{u \in
U\,|\, \tau u=\xi^\epsilon u\}$, $\epsilon=0,1,2$, where
$\xi=\exp(2\pi\sqrt{-1}/3)$. We usually denote the subspace $U(0)$
of fixed points by $U^\tau$ also.

We are interested in the subalgebra $M^\tau$. The weight $2$
subspace of $M^\tau$ is spanned by $\om$. In fact, $\om$ is the
Virasoro element of $M$ with central charge $6/5$. This means that
the subalgebra $\Vir(\om)$ generated by $\om$ is isomorphic to
$L(6/5,0)$. Note that $M$ and $M^\tau$ are completely reducible as
modules for $\Vir(\om)$, since $V_L$ possesses an invariant
positive definite hermitian form (see Subsection 5.3). Every
irreducible direct summand in $M$ or $M^\tau$ is isomorphic to
$L(6/5,h)$ for some nonnegative integer $h$. Note also that
$\sigma$ leaves $M^\tau=M(0)$ invariant and interchanges $M(1)$
and $M(2)$. Since $\sigma$ fixes $\om$, $\sigma$ acts on
$\Vir(\om)$ as the identity. Thus $M(1)$ and $M(2)$ are equivalent
$\Vir(\om)$-modules.

We now count dimension of homogeneous subspaces of $M$ of small
weights. The characters of $L(1/2,h)$, $L(7/10,h)$, and $L(6/5,h)$
are well known (cf. \cite{KR,R}). Using them, we have the first
several terms of the character of $M$:
\begin{align*}
\ch M &= \ch L(\frac{1}{2},0) \ch L(\frac{7}{10},0) + \ch
L(\frac{1}{2},\frac{1}{2}) \ch L(\frac{7}{10},\frac{3}{2})\\
&= 1+3q^2+4q^3+9q^4+12q^5+22q^6+\cdots.
\end{align*}
Comparing $\ch M$ with the character of $L(6/5,h)$, we see that
\begin{equation*}
M\cong L(\frac{6}{5},0)+2L(\frac{6}{5},2)+L(\frac{6}{5},3)+
2L(\frac{6}{5},4)+L(\frac{6}{5},6)+\cdots
\end{equation*}
as $\Vir(\om)$-modules.

The vectors $u^1$ and $u^2$ of (\ref{u1u2}) are highest weight
vectors for $\Vir(\om)$ of weight $2$. Hence the
$\Vir(\om)$-submodule generated by $u^\epsilon$ in $M(\epsilon)$
is isomorphic to $L(6/5,2)$, $\epsilon=1,2$.

Next, we study the weight $3$ subspace. The weight $3$
subspace of $M$ is of dimension $4$ and so there are nontrivial
relations among $w(\al_i)_0w(\al_j)$, $i,j \in \{0,1,2\}$. For
example,
\begin{align*}
& w(\al_1)_0w(\al_2) - w(\al_2)_0w(\al_1)\\
&\qquad=w(\al_2)_0w(\al_0) - w(\al_0)_0w(\al_2)\\
&\qquad=w(\al_0)_0w(\al_1) - w(\al_1)_0w(\al_0).
\end{align*}

Set $J=w(\al_1)_0w(\al_2) - w(\al_2)_0w(\al_1)$. In terms of the
lattice vertex operator algebra $V_L$, $J$ can be written as
\begin{equation*}
\begin{split}
J &= \frac{1}{3}\Big(\al_1(-2)\big(\al_0(-1) - \al_2(-1)\big)\\
&\qquad\qquad + \al_2(-2)\big(\al_1(-1) - \al_0(-1)\big) +
\al_0(-2)\big(\al_2(-1) - \al_1(-1)\big)\Big)\\
& \quad + \sqrt{2}\Big(\big(\al_0(-1) - \al_2(-1)\big)y(\al_1)\\
&\qquad\qquad +\big(\al_1(-1) - \al_0(-1)\big)y(\al_2) +
\big(\al_2(-1) - \al_1(-1)\big)y(\al_0)\Big).
\end{split}
\end{equation*}

Note that $(u^1)_1u^2-(u^2)_1u^1=3\sqrt{-3}J$. Note also that
$\tau J=J$, $\sigma J=-J$ and $\theta J=J$. The weight $3$
subspace of $M^\tau$ is of dimension $2$ and it is spanned by
$\om_0\om$ and $J$. Furthermore, we have $\om_1J=3J$ and
$\om_nJ=0$ for $n \ge 2$. Hence

\begin{lem} $J$ is a highest weight vector for $\Vir(\om)$ of
highest weight $3$ in $M^\tau$.
\end{lem}

The weight $4$ subspace of $M$ is of dimension $9$. By a direct
calculation, we can verify that $w(\al_i)_{-1}w(\al_j)$, $0 \le
i,j \le 2$ are linearly independent. Hence
$w(\al_i)_{-1}w(\al_j)$'s form a basis of the weight $4$ subspace
of $M$. From this it follows that the weight $4$ subspace of
$M^\tau$ is of dimension $3$. Since the weight $4$ subspace of
$\Vir(\om)\cong L(6/5,0)$ is of dimension $2$ and since the weight
$4$ subspace of the $\Vir(\om)$-submodule generated by $J$, which
is isomorphic to $L(6/5,3)$, is of dimension $1$, we conclude that
there is no highest weight vector for $\Vir(\om)$ in the weight
$4$ subspace of $M^\tau$. We have shown that

\begin{lem} $(1)$ $\{ w(\al_i)_{-1}w(\al_j)\,|\,0 \le i,j \le 2\}$
is a basis of the weight $4$ subspace of $M$.

$(2)$ There is no highest weight vector for $\Vir(\om)$ of weight
$4$ in $M^\tau$.
\end{lem}

By the above argument, we know all the irreducible direct summands
$L(6/5,h)$ with $h \le 6$ in the decomposition of $M(\epsilon)$
into a direct sum of irreducible $\Vir(\om)$-modules. Namely,
\begin{align*}
M^\tau &=M(0) \cong
L(\frac{6}{5},0)+L(\frac{6}{5},3)+L(\frac{6}{5},6)+\cdots,\\
M(1)&\cong M(2)\cong L(\frac{6}{5},2)+L(\frac{6}{5},4)+\cdots.
\end{align*}

We now consider the vertex operator algebra $\W$ generated
by $\om$ and $J$ in $M^\tau.$ Of course
$\W$ is a subalgebra of $M^{\tau}.$ We will show that $\W$ is,
in fact, equal to $M^\tau.$
The basic data are the commutation
relations of the component operators $\om_m$ and $J_n$. For the
determination of the commutation relation $[J_m,J_n]$, it is
sufficient to express $J_n J$, $0 \le n \le 5$, by using $\om$.
First of all we note that the weight $\wt J_n J = 5-n$ is at most
$5$ for $0 \le n \le 5$. Thus $J_n J$ is contained in $L(6/5, 0) +
L(6/5, 3)$, where $L(6/5, 0)$ and $L(6/5,3)$ stand for $\Vir(\om)$
and the $\Vir(\om)$-submodule generated by $J$ respectively. Since
$\sigma$ fixes every element in $\Vir(\om)$ and $\sigma J = -J$,
$\sigma$ acts as $-1$ on the $\Vir(\om)$-submodule generated by
$J$. Hence $J_n J$ is in fact contained in $\Vir(\om)$.

By a direct calculation, we have
\begin{equation}\label{JnJ}
\begin{split}
J_5 J &= -84\1,\\
J_4 J &= 0,\\
J_3 J &= -420\omega,\\
J_2 J &= -210\omega_0\omega,\\
J_1 J &= 9\omega_0\omega_0\omega - 240\omega_{-1}\omega,\\
J_0 J &= 22\omega_0\omega_0\omega_0\omega -
120\omega_0\omega_{-1}\omega.
\end{split}
\end{equation}

Note that $\{\om_0\om,\,J\}$, $\{\om_0\om_0\om,\,
\om_{-1}\om,\,\om_0 J\}$, and $\{\om_0\om_0\om_0\om,\,
\om_0\om_{-1}\om,\, \om_0\om_0 J,\, \om_{-1}J\}$ are bases of
weight $3,\,4$, and $5$ subspaces of $M^{\tau}$ respectively.

In terms of the lattice vertex operator algebra $V_L$, the vectors
$J_2 J, J_1J$, and $J_0 J$ can be written as follows.
\begin{equation*}
J_2 J = -42\sum_{i=0}^2 \al_i(-2)\al_i(-1) + 42\sqrt{2}
\sum_{i=0}^2 \al_i(-1)y(\al_i),
\end{equation*}
\begin{equation*}
\begin{split}
J_1 J &= -38 \sum_{i=0}^2 \al_i(-3)\al_i(-1) -3 \sum_{i=0}^2
\al_i(-2)^2 -8 \sum_{i=0}^2 \al_i(-1)^4\\
& \quad +6 \sum_{i=0}^2 \al_i(-1)^2 x(\al_i) + 51\sqrt{2}
\sum_{i=0}^2 \al_i(-2) y(\al_i),
\end{split}
\end{equation*}
\begin{equation*}
\begin{split}
J_0 J &= -36 \sum_{i=0}^2 \al_i(-4)\al_i(-1) - 4 \sum_{i=0}^2
\al_i(-3)\al_i(-2) - 16 \sum_{i=0}^2 \al_i(-2)\al_i(-1)^3\\
& \quad + 36 \sum_{i=0}^2 \al_i(-2)\al_i(-1)x(\al_i) +
\sum_{i=0}^2 \big(44\sqrt{2}\al_i(-3) -
4\sqrt{2}\al_i(-1)^3\big)y(\al_i).
\end{split}
\end{equation*}

We need some formulas for vertex operator algebras (cf.
\cite{FLM}), namely,
\begin{align}
[u_m,v_n] &= \sum_{k=0}^{\infty} \binom{m}{k}(u_k v)_{m+n-k},\label{commu}\\
(u_m v)_n &= \sum_{k=0}^{\infty} (-1)^k \binom{m}{k}
\big(u_{m-k}v_{n+k} - (-1)^m v_{m+n-k}u_k\big),\label{assoc}\\
(\om_0 v)_n &= -nv_{n-1}.\label{diff}
\end{align}

Using them we can obtain the commutation relations of the
component operators $\om_m$ and $J_n$.

\begin{thm} Let $L(n)=\om_{n+1}$ and $J(n)=J_{n+2}$ for $n \in \Z$,
so that the weights of these operators
are $\wt L(n)= \wt J(n)=-n$.
The
\begin{equation}
[L(m),\,L(n)]=(m-n)L(m+n)+\frac{m^3-m}{12}\cdot\frac{6}{5}
\cdot\delta_{m+n,0},\label{LL}
\end{equation}
\begin{equation}
[L(m),\,J(n)]=(2m-n)J(m+n),\label{LJ}
\end{equation}
\begin{equation}
\begin{split}
[J(m),\,J(n)]
&=(m-n)\Big(22(m+n+2)(m+n+3) + 35(m+2)(n+2)\Big)L(m+n)\\
&\qquad -120(m-n)\Big( \sum_{k \le -2} L(k)L(m+n-k) +
\sum_{k \ge -1}L(m+n-k)L(k) \Big)\\
&\qquad -\frac{7}{10}m(m^2-1)(m^2-4)\delta_{m+n,0}.
\end{split}
\label{JJ}
\end{equation}
\end{thm}

\prf The first equation holds since $\om$ is the Virasoro element
of central charge $6/5$. We know that $\om_1 J=3J$ and $\om_n J=0$
for $n \ge 2$. Hence the second equation holds. Now
\begin{align*}
(\om_{-1}\om)_{n+3}
&= \sum_{k=0}^{\infty} (-1)^k\binom{-1}{k}\Big(
\om_{-1-k}\om_{n+3+k} - (-1)^{-1}\om_{n+2-k}\om_k\Big)\\
&=\sum_{k=0}^{\infty} \Big( L(-k-2)L(n+k+2) + L(n+1-k)L(k-1)\Big)\\
&=\sum_{k \le -2} L(k)L(n-k) + \sum_{k \ge -1}L(n-k)L(k) \Big).
\end{align*}

Thus the last equation follows from (\ref{JnJ}). \quad \qed

\begin{rmk} Let $L_n=L(n)$ and $W_n=\sqrt{-1/210}J(n)$. Then
the commutation relations in the above theorem coincide with the
commutation relations (2.1) and (2.2) of \cite{BMP}. Thus $\W$ is
a $W_3$ algebra of central charge $6/5$.
\end{rmk}

Let $\lambda(m)=i(i+1)$ if $m=2i+1$ is odd and $\lambda(m)=i^2$ if
$m=2i$ is even. Let $:L(n_1)L(n_2):$ be the normal ordered
product, so that it is equal to $L(n_1)L(n_2)$ if $n_1 \le n_2$
and $L(n_2)L(n_1)$ if $n_1 \ge n_2$. Then we have another
expression of $(\om_{-1}\om)_{n+3}$. That is (cf. \cite{FZ}),
\begin{equation*}
(\om_{-1}\om)_{n+3}=\lambda(n+3)L(n) + \sum_{k \in \Z}
:L(k)L(n-k):.
\end{equation*}

\section{$20$ irreducible modules for $M^\tau$}

In this section we construct $20$ irreducible modules for
$M^\tau$. Furthermore, we calculate the action of the weight
preserving component operators $L(0)=\om_1$ and $J(0)=J_2$ on the
top levels of those irreducible modules for $M^\tau$. Recall that
$M$ has exactly $8$ inequivalent irreducible modules $M_k^i$,
$W_k^i$, $i=0,a,b,c$. Let $(U, Y_U)$ be one of those irreducible
$M$-modules. Following \cite{DLM}, we consider a new $M$-module
$(U\circ\tau, Y_{U\circ\tau})$ such that $U\circ\tau = U$ as
vector spaces and
\begin{equation*}
Y_{U\circ\tau}(v,z)=Y_U(\tau v,z) \quad \mbox{for } v \in M.
\end{equation*}
Then $U \mapsto U\circ\tau$ induces a permutation on the set of
irreducible $M$-modules. If $U$ and $U\circ\tau$ are equivalent
$M$-modules, $U$ is said to be $\tau$-stable. By the definition,
we have $U\circ\tau^2=(U\circ\tau)\circ\tau$. The following lemma
is a straightforward consequence of the definition of $M_k^i$ and
$W_k^i$.

\begin{lem} $(1)$ $M_k^0\circ\tau=M_k^0$ and
$W_k^0\circ\tau=W_k^0$.

$(2)$ $M_k^a\circ\tau=M_k^c$, $M_k^c\circ\tau=M_k^b$, and
$M_k^b\circ\tau=M_k^a$.

$(3)$ $W_k^a\circ\tau=W_k^c$, $W_k^c\circ\tau=W_k^b$, and
$W_k^b\circ\tau=W_k^a$.
\end{lem}

Here $W_k^0\circ\tau=W_k^0$ means that there exists a linear
isomorphism $\phi(\tau) : W_k^0 \longrightarrow W_k^0$ such that
$\phi(\tau) Y_{W_k^0}(v,z)\phi(\tau)^{-1} = Y_{W_k^0}(\tau v,z)$
for all $v \in M$. The automorphism $\tau$ of $V_L$ fixes
$\omega^3$ and so $W_k^0$ is invariant under $\tau$. Hence we can
take $\tau$ as $\phi(\tau)$. Note also that $\tau
Y(v,z)\tau^{-1}=Y(\tau v,z)$ for all $v \in M=M_k^0$ since $\tau
\in \Aut (M)$.

\subsection{Irreducible $M^\tau$-modules in untwisted $M$-modules}

We first find $8$ irreducible $M^\tau$-modules inside the $8$
irreducible modules for $M$. Recall that $M(\epsilon)=\{ v \in
M_k^0\,|\, \tau v= \xi^\epsilon v\}$. Likewise, set
$W(\epsilon)=\{ v \in W_k^0\,|\, \tau v= \xi^\epsilon v\}$. From
Lemma 4.1, \cite[Theorem 4.4]{DM}, and \cite[Theorem 6.14]{DY}, we
see that $M(\epsilon)$ and $W(\epsilon)$ are inequivalent
irreducible $M^\tau$-modules for $\epsilon=0,1,2$. Note that
$M_k^i$, $i=a,b,c$ are equivalent irreducible $M^\tau$-modules and
that $W_k^i$, $i=a,b,c$ are also equivalent irreducible
$M^\tau$-modules by \cite[Theorem 6.14]{DY}. Hence we obtain 8
inequivalent irreducible $M^\tau$-modules.

The top levels, that is, the weight subspaces of the smallest
weights of $M(0)$, $M(1)$, and $M(2)$ are $\C\1$, $\C u^1$, and
$\C u^2$ respectively. The top levels of $W(0)$, $W(1)$, and
$W(2)$ are
\begin{equation*}
\C(y(\al_1)+y(\al_2)+y(\al_0)), \quad \C(\al_1(-1)-\xi\al_2(-1)),
\quad and \quad \C(\al_1(-1)-\xi^2\al_2(-1))
\end{equation*}
respectively. Moreover, the top levels of $M_k^c$ and $W_k^c$ are
\begin{equation*}
\C(e^{\beta_1/2} - e^{-\beta_1/2}) \quad \mbox{and} \quad
\C(e^{\beta_1/2} + e^{-\beta_1/2})
\end{equation*}
respectively. All of those top levels are of dimension one.

Next, we deal with the action of $L(0)$ and $J(0)$ on those top
levels. The operator $L(0)$ acts as multiplication by the weight
of each top level. For the calculation of the action of $J(0)$, we
first notice that
\begin{equation*}
[w(\al_i)_1, w(\al_j)_1] = (w(\al_i)_0 w(\al_j))_2 + (w(\al_i)_1
w(\al_j))_1
\end{equation*}
by (\ref{commu}). Since $w(\al_i)_1 w(\al_j) = w(\al_j)_1
w(\al_i)$, it follows that
\begin{align*}
J(0) &=(w(\al_1)_0 w(\al_2))_2 - (w(\al_2)_0 w(\al_1))_2\\
&=[w(\al_1)_1, w(\al_2)_1] - [w(\al_2)_1, w(\al_1)_1].
\end{align*}

Using this formula it is relatively easy to calculate the
eigenvalue for the action of $J(0)$ on each of the $8$ top levels.
The results are collected in Table \ref{top1}.

\begin{table}[h]
\begin{center}
\caption{irreducible $M^\tau$-modules in $M_k^i$ and $W_k^i$}
\label{top1}
\begin{tabular}{l|l|l|l}
irred. module & top level & $L(0)$ & $J(0)$\\
\hline
$M(0)$ & $\C\1$ & $0$ & $0$\\
$M(1)$ & $\C u^1$ & $2$ & $-12\sqrt{-3}$\\
$M(2)$ & $\C u^2$ & $2$ & $12\sqrt{-3}$\\
$W(0)$ & $\C(y(\al_1)+y(\al_2)+y(\al_0))$ & $\frac{8}{5}$
& $0$\\
$W(1)$ & $\C(\al_1(-1) - \xi\al_2(-1))$ & $\frac{3}{5}$
& $2\sqrt{-3}$\\
$W(2)$ & $\C(\al_1(-1) - \xi^2\al_2(-1))$ & $\frac{3}{5}$
& $-2\sqrt{-3}$\\
$M_k^c$ & $\C(e^{\be_1/2}-e^{-\be_1/2})$ & $\frac{1}{2}$ & $0$
\\
$W_k^c$ & $\C(e^{\be_1/2}+e^{-\be_1/2})$ & $\frac{1}{10}$ & $0$
\end{tabular}
\end{center}
\end{table}

\subsection{Irreducible $M^\tau$-modules in $\tau$-twisted
$M$-modules}

Using \cite{DLM}, we show that there are exactly two inequivalent
irreducible $\tau$-twisted (resp. $\tau^2$-twisted) $M$-modules.
Moreover, we find $3$ inequivalent irreducible $M^\tau$-modules in
each of the irreducible $\tau$-twisted (resp. $\tau^2$-twisted)
$M$-modules. Those irreducible $\tau$-twisted (resp.
$\tau^2$-twisted) $M$-modules will in turn be constructed inside
irreducible $\tau$-twisted (resp. $\tau^2$-twisted) $V_L$-modules.
Basic references to twisted modules for lattice vertex operator
algebras are \cite{D2, DL, L}. The argument here is similar to
that in \cite[Section 6]{KLY2}.

We follow \cite{DL} with $L=\sqrt{2}A_2$, $p=3$, $q=6$, and
$\nu=\tau$. Let $\h = \C\otimes_{\Z}L$ and extend the
$\Z$-bilinear form $\la\cdot,\cdot\ra$ on $L$ to $\h$ linearly.
Set
\begin{equation*}
h_1 = \frac{1}{3}(\be_1 + \xi^2 \be_2 + \xi \be_0), \qquad h_2=
\frac{1}{3}(\be_1 + \xi \be_2 + \xi^2 \be_0).
\end{equation*}
Then $\tau h_j = \xi^j h_j$, $\la h_1, h_1 \ra = \la h_2, h_2 \ra
= 0$, and $\la h_1, h_2 \ra=2$. Moreover,
$\be_i=\xi^{i-1}h_1+\xi^{2(i-1)}h_2$, $i=0,1,2$. For $n \in \Z$,
set
\begin{equation*}
\h_{(n)} = \{\al \in \h\, |\,\tau\al=\xi^n\al\}.
\end{equation*}

Since $\tau$ is fixed-point-free on $L$, it follows that
$\h_{(0)}=0$. Furthermore, $\h_{(1)}=\C h_1$ and $\h_{(2)}=\C
h_2$. For $\al \in \h$, we denote by $\al_{(n)}$ the component of
$\al$ in $\h_{(n)}$. Thus $(\be_i)_{(1)}=\xi^{i-1}h_1$ and
$(\be_i)_{(2)}=\xi^{2(i-1)}h_2$ for $i=0,1,2$.

Define the $\tau$-twisted affine Lie algebra to be
\begin{equation*}
\hat{\h}[\tau]=\Big(\bigoplus_{n\in \Z}\h_{(n)}\otimes t^{n/3}\Big)
\oplus \C c
\end{equation*}
with the bracket
\begin{equation*}
[x\otimes t^m, y\otimes t^n]= m \la x,y\ra \delta_{m+n,0}c
\end{equation*}
for $x\in \h_{(3m)}$, $y\in \h_{(3n)}$, $m,n\in (1/3)\Z$, and
$[c,\hat{\h}[g]]=0$. The isometry $\tau$ acts on $\hat{\h}[\tau]$
by $\tau (x \otimes t^m)=\xi^m x\otimes t^m$ and $\tau(c)=c$. Set
\begin{equation*}
\hat{\h}[\tau]^+=\bigoplus_{n>0}\h_{(n)}\otimes t^{n/3},\quad
\hat{\h}[\tau]^-=\bigoplus_{n<0}\h_{(n)}\otimes t^{n/3},\quad
\text{and} \quad \hat{\h}[\tau]^0 =\C c
\end{equation*}
and consider the $\hat{\h}[\tau]$-module
\begin{equation*}
S[\tau]=U(\hat{\h}[\tau])\otimes_{U(\hat{\h}[\tau]^+\oplus
\hat{\h}[\tau]^0)}\C
\end{equation*}
induced from the $\hat{\h}[\tau]^+\oplus
\hat{\h}[\tau]^0$-module $\C$, where $\hat{\h}[\tau]^+$ acts
trivially on $\C$ and $c$ acts as $1$ on $\C$.

We define the weight in $S[\tau]$ by
\begin{equation*}
\wt (x\otimes t^n)=-n \quad \text{and}\quad \wt 1 =\frac{1}{9},
\end{equation*}
where $n \in (1/3)\Z$ and $x \in \h_{(3n)}$ (cf. \cite[(4.6),
(4.10)]{DL}). By the weight gradation $S[\tau]$ becomes a
$(1/3)\Z$-graded space. Its character is
\begin{equation}\label{chs}
\ch S[\tau] = q^{1/9}\prod_{n=1}^{\infty} (1-q^n)/
\prod_{n=1}^{\infty} (1-q^{n/3}).
\end{equation}

For $\al \in \h$ and $n \in (1/3)\Z$, denote by $\al(n)$ the
operator on $S[\tau]$ induced by $\al_{(3n)}\otimes t^n$. Then, as
a vector space $S[\tau]$ can be identified with a polynomial
algebra with variables $h_1(1/3+n)$ and $h_2(2/3+n)$, $n \in \Z$.
The weight of the operator $h_j(j/3+n)$ is $-j/3-n$.

The alternating $\Z$-bilinear map $c_0^\tau : L \times L \to
\Z/6\Z$ defined by \cite[(2.10)]{DL} is such that
\begin{equation*}
c_0^\tau(\al,\beta) = \sum_{r=0}^2 (3+2r)\la\tau^r\al, \beta\ra +
6\Z.
\end{equation*}

In our case $\sum_{r=0}^2 \tau^r\al = 0$, since $\tau$ is
fixed-point-free on $L$. Moreover, we can verify that
\begin{equation*}
\sum_{r=0}^2 r\la\tau^r \beta_i, \beta_j\ra =
\begin{cases}
\pm 6 & \text{if\ } \tau\beta_i \ne \beta_j\\
0 & \text{if\ } \tau\beta_i = \beta_j.
\end{cases}
\end{equation*}
Hence $c_0^\tau(\al,\beta) = 0$ for all $\al, \beta \in L$.
This means that
the central extension
\begin{equation}\label{ext2}
1 \longrightarrow \la\kappa_6\ra \longrightarrow \hat{L}_\tau
\overset{-}{\longrightarrow} L \longrightarrow 1
\end{equation}
determined by the commutator condition
$aba^{-1}b^{-1}=\kappa_6^{c_0^\tau(\bar{a},\bar{b})}$ splits.

We consider the relation between two central extensions $\hat{L}$
of (\ref{ext1}) and $\hat{L}_\tau$ of (\ref{ext2}). Since both of
$\hat{L}$ and $\hat{L}_\tau$ are split extensions, we use the same
symbol $e^\al$ to denote both of an element in $\hat{L}$ and an
element in $\hat{L}_\tau$ which correspond naturally to $\al \in
L$. Actually, in Section 2 we choose $e^\al \in \hat{L}$ so that
the multiplication in $\hat{L}$ is $e^\al \times e^\be =
e^{\al+\be}$. Also we can choose $e^\al \in \hat{L}_\tau$ such
that the multiplication $e^\al \times_\tau e^\be$ in
$\hat{L}_\tau$ is related to the multiplication in $\hat{L}$ by
(cf. \cite[(2.4)]{DL})
\begin{equation}\label{mult}
e^\al \times e^\be = \kappa_6^{\ve_0(\al,\be)} e^\al \times_\tau
e^\be,
\end{equation}
where the $\Z$-linear map $\varepsilon_0 : L \times L \to \Z/6\Z$
is defined by \cite[(2.13)]{DL}. In our case
\begin{equation}\label{ve0}
\ve_0(\al,\be)=-\la\tau^{-1}\al,\be\ra + 6\Z.
\end{equation}

As in Section 2, we usually write $e^\al e^\be=e^{\al+\be}$ to
denote the product
of $e^\al$ and $e^\be$ in $\hat{L}$. Note, for example, that
the inverse of $e^{\be_1}$ in $\hat{L}$ is $e^{-\be_1}$, while
the inverse of $e^{\be_1}$ in $\hat{L}_\tau$ is $\kappa_3^2
e^{-\be_1}$.

The automorphism $\tau$ of $L$ lifts to an automorphism
$\hat{\tau}$ of $\hat{L}$ such that $\hat{\tau}(e^\al) =
e^{\tau\al}$ and $\hat{\tau}(\kappa_6)=\kappa_6$. Since $\ve_0$ is
$\tau$-invariant, we can also think $\hat{\tau}$ to be an automorphism
of $\hat{L}_\tau$ in a similar way. By abuse of notation we shall
denote $\hat{\tau}$ by simply $\tau$ also.

We have $(1-\tau)L=\spn_{\Z}\{ \be_1-\be_2, \be_1+2\be_2\}$. The
quotient group $L/(1-\tau)L$ is of order $3$ and generated by
$\be_1+(1-\tau)L$. Now $K=\{
a^{-1}\tau(a)\,|\,a\in\hat{L}_\tau\}$ is a central subgroup
of $\hat{L}_\tau$ with $\overline{K}= (1-\tau)L$ and $K \cap
\la\kappa_6\ra=1$. Here note that $a^{-1}$ is the inverse of $a$
in $\hat{L}_\tau$ and $a^{-1}\tau(a)$ is the product $a^{-1}
\times_\tau \tau(a)$ in $\hat{L}_\tau$. In $\hat{L}_\tau$
we can verify that
\begin{equation*}
e^{3\be_1}=(e^{\be_0-\be_1})^{-1} \times_\tau
\tau(e^{\be_0-\be_1}) \in K.
\end{equation*}
Since
\begin{equation*}
\kappa_3 e^{\be_1} \times_\tau \kappa_3 e^{\be_1} \times_\tau
\kappa_3 e^{\be_1} = e^{3\be_1} \qquad \text{and}\qquad
\kappa_3 e^{\be_1} \times_\tau \kappa_3 e^{-\be_1}=1,
\end{equation*}
it follows that
\begin{equation*}
\hat{L}_\tau/K=\{ K, \kappa_3 e^{\be_1}K, \kappa_3 e^{-\be_1}K\}
\times \la\kappa_6\ra K/K \cong \Z_3 \times \Z_6.
\end{equation*}
For $j=0,1,2$, define a linear character $\chi_j : \hat{L}_\tau/K
\to \C^{\times}$ by
\begin{equation*}
\chi_j(\kappa_6)=\xi_6, \quad \chi_j(\kappa_3 e^{\be_1}K)=\xi^j,
\quad \text{and} \quad \chi_j(\kappa_3 e^{-\be_1}K)=\xi^{-j},
\end{equation*}
where $\xi_6=\exp(2\pi\sqrt{-1}/6)$. Let $T_{\chi_j}$ be the one
dimensional $\hat{L}_\tau/K$-module affording the character
$\chi_j$. As an $\hat{L}_\tau$-module, $K$ acts trivially on
$T_{\chi_j}$. Since $\sum_{r=0}^2 \tau^r \al=0$ for $\al \in L$,
those $T_{\chi_j}$, $j=0,1,2$, are the irreducible
$\hat{L}_\tau$-modules constructed in \cite[Section 6]{L}.

Let
\begin{equation*}
V_L^{T_{\chi_j}}=V_L^{T_{\chi_j}}(\tau)=S[\tau]\otimes T_{\chi_j}
\end{equation*}
and define the $\tau$-twisted vertex operator $Y^\tau(\cdot,
z): V_L \to \End(V_L^{T_{\chi_j}})\{z\}$ as in \cite{DL}. For
$a\in\hat{L}$, define
\begin{equation*}
Y^\tau(a,z)=3^{-\la\bar{a},\bar{a}\ra/2} \phi(\bar{a})
E^-(-\bar{a},z)E^{+}(-\bar{a},z) a z^{-\la\bar{a},\bar{a}\ra/2},
\end{equation*}
where
\begin{equation}
E^{\pm}(\al,z)=\exp\left(\sum_{n\in (1/3)\mathbb{Z}_{\pm}}
\frac{\al(n)}{n}z^{-n}\right),
\end{equation}
\begin{equation}
\phi(\alpha)= (1-\xi^2)^{\langle \tau\alpha,\alpha\rangle},
\end{equation}
and $a\in \hat{L}$ acts on $T_{\chi_j}$ through the set theoretic
identification between $\hat{L}$ and $\hat{L}_\tau$. Here we
denote $\sigma(\al)$ of \cite[(4.35)]{DL} by $\phi(\al)$. For
$v=\alpha^1(-n_1)\cdots\alpha^k(-n_k)\cdot\iota(a)\in V_L$ with
$\alpha^1,\ldots,\alpha^k \in\mathfrak{h}$ and $n_1,\ldots,n_k \in
\Z_{> 0}$, set
\begin{equation*}
W(v,z)=\nop
\left(\frac{1}{(n_1-1)!}\left(\frac{d}{dz}\right)^{n_1-1}
\alpha^1(z)\right)\cdot\cdot\cdot\left(\frac{1}{(n_k-1)!}
\left(\frac{d}{dz}\right)^{n_k-1} \alpha^k(z)\right)Y^{\tau}(a,z)
\nop,
\end{equation*}
where $\al(z)=\sum_{n\in (1/3)\Z} \al(n)z^{-n-1}$.
Define constants $c_{mn}^i\in\C$ for $m, n\ge 0$ and
$i=0,1,2$ by
\begin{gather*}
\sum_{m,n\ge 0}c_{mn}^0x^my^n=-\frac{1}{2}\sum_{r=1}^{2}{\rm
log}\left(\frac
{(1+x)^{1/3}-\xi^{-r}(1+y)^{1/3}}{1-\xi^{-r}}\right),\\
\sum_{m,n\ge 0}c_{mn}^ix^my^n=\frac{1}{2}{\rm log}\left( \frac
{(1+x)^{1/3}-\xi^{-i}(1+y)^{1/3}}{1-\xi^{-i}}\right)\ \text{ for}\
\ i\ne0.
\end{gather*}

Let $\{\gm_1,\gm_2\}$ be an orthonormal basis of
$\mathfrak{h}$ and set
\begin{equation*}
\Delta_z=\sum_{m,n\ge 0}\displaystyle{\sum^{2}_{i=0}}\
\displaystyle{ \sum^2_{j=1}}
c_{mn}^i(\tau^{-i}\gm_j)(m)\gm_j(n)z^{-m-n}.
\end{equation*}
Then for $v\in V_L,$  $Y^{\tau}(v,z)$ is defined by
\begin{equation*}
Y^{\tau}(v,z)=W(e^{\Delta_z}v,z).
\end{equation*}

We extend the action of $\tau$ to $V_L^{T_{\chi_j}}$ so that
$\tau$ is the identity on $T_{\chi_j}$. The weight of every
element in $T_{\chi_j}$ is defined to be $0$. Then the character
of $V_L^{T_{\chi_j}}$ is identical with that of $S[\tau]$.

By \cite[Theorem 7.1]{DL}, $(V_L^{T_{\chi_j}}(\tau),
Y^\tau(\cdot,z))$, $j=0,1,2$ are inequivalent irreducible
$\tau$-twisted $V_L$-modules. Now among the $12$ irreducible
$V_L$-modules $V_{L^{(i,j)}}$, $i\in \{0,a,b,c\}$ and $j \in
\{0,1,2\}$, the $\tau$-stable irreducible modules are
$V_{L^{(0,j)}}$, $j \in \{0,1,2\}$. Hence by \cite[Theorem
10.2]{DLM}, we conclude that $(V_L^{T_{\chi_j}}(\tau),
Y^\tau(\cdot,z))$, $j=0,1,2$, are all the inequivalent irreducible
$\tau$-twisted $V_L$-modules. The isometry $\theta$ of $(L,
\la\cdot,\cdot\ra)$ induces a permutation on
$V_L^{T_{\chi_j}}(\tau)$, $j=0,1,2$. In fact, the permutation
leaves $V_L^{T_{\chi_0}}(\tau)$ invariant and interchanges
$V_L^{T_{\chi_1}}(\tau)$ and $V_L^{T_{\chi_2}}(\tau)$.

Since $M^\tau \otimes M_t^0$ is contained in the subalgebra
$(V_L)^\tau$ of fixed points of $\tau$ in $V_L$, we can deal with
$(V_L^{T_{\chi_j}}(\tau),Y^\tau(\cdot,z))$ as an $M^\tau \otimes
M_t^0$-module. We will find $6$ irreducible $M^\tau$-modules
inside $V_L^{T_{\chi_j}}(\tau)$. Recall that $\om$, $\om^3$, and
$\widetilde{\om}=\om+\om^3$ are the Virasoro element of $M^\tau$,
$M_t^0$, and $V_L$ respectively. Our main tool is a careful study
of the action of $\om_1$ on homogeneous subspaces of
$V_L^{T_{\chi_j}}(\tau)$ of small weights. Here we denote by $u_n$
the coefficient of $z^{-n-1}$ in the twisted vertex operator
$Y^\tau(u,z)=\sum u_n z^{-n-1}$ associated with a vector $u$ in
$V_L$. The weight in $V_L^{T_{\chi_j}}(\tau)$ defined above is
exactly the eigenvalue for $\widetilde{\om}_1$ (cf. \cite[(6.10),
(6.28)]{DL}).

The character of $V_L^{T_{\chi_j}}(\tau)$ is equal to the
character of $S[\tau]$ (cf. (\ref{chs})). Its first several terms
are
\begin{equation*}
\ch V_L^{T_{\chi_j}}(\tau) = \ch S[\tau] = q^{1/9} + q^{1/9+1/3} +
2q^{1/9+2/3} + 2q^{1/9+1} + 4q^{1/9+4/3} + \cdots.
\end{equation*}

Fix a nonzero vector $v \in T_{\chi_j}$. We can choose a basis
of each homogeneous subspace of $V_L^{T_{\chi_j}}(\tau)$ of
weight at most $1/9+4/3$ as in Table \ref{basis}.

\begin{table}[h]
\begin{center}
\caption{Basis of homogeneous subspace in $V_L^{T_{\chi_j}}(\tau)$}
\label{basis}
\begin{tabular}{l|l}
weight & \hspace{4cm} basis\\
\hline
$\frac{1}{9}$ & $1 \otimes v$\\
$\frac{1}{9}+\frac{1}{3}$ & $h_2(-\frac{1}{3}) \ot v$\\
$\frac{1}{9}+\frac{2}{3}$ & $h_1(-\frac{2}{3}) \ot v$, \quad
$h_2(-\frac{1}{3})^2 \ot v$\\
$\frac{1}{9}+1$ & $h_1(-\frac{2}{3})h_2(-\frac{1}{3}) \ot v$, \quad
$h_2(-\frac{1}{3})^3 \ot v$\\
$\frac{1}{9}+\frac{4}{3}$ & $h_2(-\frac{4}{3}) \ot v$, \quad
$h_1(-\frac{2}{3})^2 \ot v$, \quad
$h_1(-\frac{2}{3})h_2(-\frac{1}{3})^2 \ot v$, \quad
$h_2(-\frac{1}{3})^4 \ot v$
\end{tabular}
\end{center}
\end{table}

We need to know the action of $\om_1$ on those bases. For this
purpose, notice that
\begin{equation*}
Y^\tau(e^{\pm\beta_i}, z)=-\frac{1}{27}E^-(\mp\beta_i, z)
E^+(\mp\beta_i, z)\xi^{\pm j}z^{-2}, \qquad i=0,1,2,
\end{equation*}
since $\phi(\pm\beta_i)=-\xi/3$ and since $e^{\pm\beta_i}$ acts on
$T_{\chi_j}$ as a multiplication by $\chi_j(e^{\pm\beta_i})
=\xi^{\pm j - 1}$ for $i,j=0,1,2$. The image of the vectors in
Table \ref{basis} under the operator $\om_1$ are calculated as
follows.

\begin{equation*}
\om_1 (1\otimes v) =\big(
\frac{1}{15}+\frac{1}{45}(\xi^j+\xi^{-j})\big) 1\otimes v,
\end{equation*}
\begin{equation*}
\om_1 (h_2(-\frac{1}{3})\otimes v) = \big(
\frac{4}{15}-\frac{1}{9}(\xi^j+\xi^{-j})\big) h_2(-\frac{1}{3})
\otimes v, \end{equation*}
\begin{equation*} \om_1
(h_1(-\frac{2}{3})\otimes v) =\big(
\frac{7}{15}-\frac{2}{45}(\xi^j+\xi^{-j})\big)
h_1(-\frac{2}{3})\otimes v - \frac{1}{5}(\xi^j-\xi^{-j})
h_2(-\frac{1}{3})^2\otimes v,
\end{equation*}
\begin{equation*}
\om_1 (h_2(-\frac{1}{3})^2 \otimes v) =\big(
\frac{7}{15}+\frac{7}{45}(\xi^j+\xi^{-j})\big)
h_2(-\frac{1}{3})^2\otimes v + \frac{2}{15}(\xi^j-\xi^{-j})
h_1(-\frac{2}{3})^2\otimes v,
\end{equation*}
\begin{equation*}
\om_1 (h_1(-\frac{2}{3})h_2(-\frac{1}{3})\otimes v) =\big(
\frac{2}{3}+\frac{2}{9}(\xi^j+\xi^{-j})\big)
h_1(-\frac{2}{3})h_2(-\frac{1}{3})\otimes v +
\frac{1}{5}(\xi^j-\xi^{-j}) h_2(-\frac{1}{3})^3\otimes v,
\end{equation*}
\begin{equation*}
\om_1 (h_2(-\frac{1}{3})^3 \otimes v) =\big(
\frac{2}{3}+\frac{1}{45}(\xi^j+\xi^{-j})\big)
h_2(-\frac{1}{3})^3\otimes v - \frac{2}{5}(\xi^j-\xi^{-j})
h_1(-\frac{2}{3})h_2(-\frac{1}{3})\otimes v,
\end{equation*}
\begin{align*}
\om_1 (h_2(-\frac{4}{3}) \otimes v)
&=\frac{13}{15}h_2(-\frac{4}{3})\otimes v\\
&\qquad+(\xi^j+\xi^{-j})\big( -\frac{1}{90}h_2(-\frac{4}{3}) -
\frac{3}{10}h_1(-\frac{2}{3})h_2(-\frac{1}{3})^2\big)\otimes v\\
&\qquad+(\xi^j-\xi^{-j})\big( -\frac{1}{20}h_1(-\frac{2}{3})^2 -
\frac{3}{20}h_2(-\frac{1}{3})^4\big)\otimes v,
\end{align*}
\begin{align*}
\om_1 (h_1(-\frac{2}{3})^2 \otimes v)
&=\frac{13}{15}h_1(-\frac{2}{3})^2\otimes v\\
&\qquad+(\xi^j+\xi^{-j})\big( -\frac{1}{90}h_1(-\frac{2}{3})^2 +
\frac{3}{10}h_2(-\frac{1}{3})^4)\big)\otimes v\\
&\qquad+(\xi^j-\xi^{-j})\big(\frac{1}{15}h_2(-\frac{4}{3}) +
\frac{1}{5}h_1(-\frac{2}{3})h_2(-\frac{1}{3})^2\big)\otimes v,
\end{align*}
\begin{align*}
\om_1 (h_1(-\frac{2}{3})h_2(-\frac{1}{3})^2 \otimes v)
&=\frac{13}{15}h_1(-\frac{2}{3})h_2(-\frac{1}{3})^2\otimes v\\
&\qquad+(\xi^j+\xi^{-j})\big( -\frac{2}{15}h_2(-\frac{4}{3}) -
\frac{14}{45}h_1(-\frac{2}{3})h_2(-\frac{1}{3})^2\big)\otimes v\\
&\qquad-\frac{1}{15}(\xi^j-\xi^{-j})h_1(-\frac{2}{3})^2\otimes v,
\end{align*}
\begin{align*}
\om_1 (h_2(-\frac{1}{3})^4 \otimes v)
&=\frac{13}{15}h_2(-\frac{1}{3})^4\otimes v\\
&\qquad+(\xi^j+\xi^{-j})\big(\frac{2}{5}h_1(-\frac{2}{3})^2 -
\frac{1}{9}h_2(-\frac{1}{3})^4\big)\otimes v\\
&\qquad+\frac{4}{15}(\xi^j-\xi^{-j})h_2(-\frac{4}{3})\otimes v,\\
\end{align*}

The decomposition of $V_L^{T_{\chi_j}}(\tau)$ as a $\tau$-twisted
$M\otimes M_t^0$-module was studied in \cite{KLY2}. The outline of
the argument is as follows. For $j=0,1,2$, the vectors
\begin{equation*}
1 \otimes v, \quad h_1(-\frac{2}{3})\otimes v +
(\xi^j-\xi^{-j})h_2(-\frac{1}{3})^2\otimes v, \quad
h_2(-\frac{1}{3})^2 \otimes v +
\frac{2}{3}(\xi^j-\xi^{-j})h_1(-\frac{2}{3})\otimes v
\end{equation*}
are simultaneous eigenvectors for $\om_1$ and $(\om^3)_1$. Denote
by $k_1$ and $k_2$ the eigenvalues for $\om_1$ and $(\om^3)_1$
respectively. Then the pairs $(k_1,k_2)$ are

\begin{center}
\begin{tabular}{c|ccc}
j & 0 & 1 & 2\\
\hline
$1\otimes v$ & $(\frac{1}{9},0)$ & $(\frac{2}{45}, \frac{1}{15})$
& $(\frac{2}{45}, \frac{1}{15})$\\
$h_1(-\frac{2}{3})\otimes v +
(\xi^j-\xi^{-j})h_2(-\frac{1}{3})^2\otimes v$ &
$(\frac{17}{45},\frac{2}{5})$ & $(\frac{1}{9},\frac{2}{3})$
& $(\frac{1}{9},\frac{2}{3})$\\
$h_2(-\frac{1}{3})^2 \otimes v +
\frac{2}{3}(\xi^j-\xi^{-j})h_1(-\frac{2}{3})\otimes v$ &
$(\frac{7}{9},0)$ & $(\frac{32}{45},\frac{1}{15})$ &
$(\frac{32}{45},\frac{1}{15})$
\end{tabular}
\end{center}

We first discuss the decomposition of $V_L^{T_{\chi_j}}(\tau)$
into a direct sum of irreducible $M_t^0$-modules. We use the
classification of irreducible $M_t^0$-modules \cite{KMY} and their
fusion rules \cite{M2}. Note also that the vector $y(\al_1) +
y(\al_2)+y(\al_0)$ in $(V_L)^\tau$ is an eigenvector for $\om_1$
of eigenvalue $8/5$. Hence $(V_L)^\tau$ contains the
$\Vir(\om)\otimes M_t^0$-submodule generated by the vector, which
is isomorphic to
\begin{equation*}
L(\frac{6}{5},\frac{8}{5})\otimes\Big(L(\frac{4}{5},\frac{2}{5}) +
L(\frac{4}{5},\frac{7}{5})\Big).
\end{equation*}

Set
\begin{align*}
M^0_T(\tau)&=\{ u\in V_L^{T_{\chi_0}}(\tau)\,|\,(\om^3)_1 u=0
\},\\
W^0_T(\tau)&=\{ u\in V_L^{T_{\chi_0}}(\tau)\,|\,(\om^3)_1
u=\frac{2}5 u \}.
\end{align*}
Moreover, for $j=1,2$ set
\begin{align*}
M^j_T(\tau)&=\{ u\in V_L^{T_{\chi_j}}(\tau)\,|\,(\om^3)_1
u=\frac{2}3 u\},\\
W^j_T(\tau)&=\{ u\in V_L^{T_{\chi_j}}(\tau)\,|\,(\om^3)_1
u=\frac{1}{15} u \}.
\end{align*}
Then, by \cite[Proposition 6.8]{KLY2}, $M_T^j(\tau)$ and
$W_T^j(\tau)$, $j=0,1,2$, are irreducible $\tau$-twisted
$M$-modules. Furthermore, for $j=0,1,2$,
\begin{equation*}
V_L^{T_{\chi_j}}(\tau) \cong M_T^j(\tau)\otimes M_t^j \oplus
W_T^j(\tau)\otimes W_t^j
\end{equation*}
as $\tau$-twisted $M \otimes M_t^0$-modules.

There are at most two inequivalent irreducible $\tau$-twisted
$M$-modules by Lemma 4.1 and \cite[Theorem 10.2]{DLM}. Then,
looking at the smallest weight of $M_T^j(\tau)$ and $W_T^j(\tau)$,
we have that $M_T^0(\tau) \cong M_T^1(\tau) \cong M_T^2(\tau)$ and
$W_T^0(\tau) \cong W_T^1(\tau) \cong W_T^2(\tau)$ and that
$M_T^0(\tau)\not\cong W_T^0(\tau)$ as $\tau$-twisted $M$-modules.
We denote $M_T^0(\tau)$ by $M_T(\tau)$ and $W_T^0(\tau)$ by
$W_T(\tau)$. We conclude that there are exactly two inequivalent
irreducible $\tau$-twisted $M$-modules, which are represented by
$M_T(\tau)$ and $W_T(\tau)$. As $\tau$-twisted $M\otimes
\Vir(\om^3)$-modules, we have
\begin{align}
V_L^{T_{\chi_0}}(\tau) & \cong M_T(\tau)\otimes
\left(L(\frac{4}5,0)+L(\frac{4}5,3)\right) \oplus W_T(\tau)\otimes
\left(L(\frac{4}5,\frac{2}{5})+L(\frac{4}5,\frac{7}{5})\right),
\label{decomp1}\\
V_L^{T_{\chi_1}}(\tau) & \cong V_L^{T_{\chi_2}}(\tau) \cong
M_T(\tau)\otimes L(\frac{4}{5},\frac{2}{3}) \oplus
W_T(\tau)\otimes L(\frac{4}{5},\frac{1}{15}). \label{decomp2}
\end{align}

The first several terms of the characters of $M_T(\tau)$ and
$W_T(\tau)$ are
\begin{align*}
\ch M_T(\tau) &=q^{\frac{1}{9}}+q^{\frac{1}{9}+\frac{2}{3}}+
q^{\frac{1}{9}+1}+q^{\frac{1}{9}+\frac{4}{3}}+\cdots,\\
\ch W_T(\tau) &=q^{\frac{2}{45}}+
q^{\frac{2}{45}+\frac{1}{3}}+q^{\frac{2}{45}+\frac{2}{3}}+
q^{\frac{2}{45}+1}+2q^{\frac{2}{45}+\frac{4}{3}}+\cdots .
\end{align*}

For $\epsilon=0,1,2$, let
\begin{align*}
M_T(\tau)(\epsilon) &=\{ u \in M_T(\tau)\,|\,\tau u=\xi^\epsilon u\},\\
W_T(\tau)(\epsilon) &=\{ u \in W_T(\tau)\,|\,\tau u=\xi^\epsilon
u\}.
\end{align*}
Those $6$ modules for $M^\tau$ are inequivalent irreducible
modules by \cite[Theorem 2]{MT}. Their top levels are of dimension
one. Those top levels and the eigenvalues for the action of
$L^\tau(0)=\om_1$ and $J^\tau(0)=J_2$ are collected in Table
\ref{toptau}.

\begin{table}[h]
\begin{center}
\caption{irreducible $M^\tau$-modules in $M_T(\tau)$ and
$W_T(\tau)$} \label{toptau}
\begin{tabular}{l|l|l|l}
irred. module & top level & $L^\tau(0)$ & $J^\tau(0)$\\
\hline $M_T(\tau)(0)$ & $\C 1\otimes v$ & $\frac{1}{9}$ &
$\frac{14}{81}\sqrt{-3}$\\
$M_T(\tau)(1)$ & $\C h_2(-\frac{1}{3})^2\otimes v$ &
$\frac{1}{9}+\frac{2}{3}$ &
$-\frac{238}{81}\sqrt{-3}$\\
$M_T(\tau)(2)$ & $\C (\frac{4}{3}h_1(-\frac{2}{3})^2\otimes v +
h_2(-\frac{1}{3})^4\otimes v)$ & $\frac{1}{9}+\frac{4}{3}$ &
$\frac{374}{81}\sqrt{-3}$\\
$W_T(\tau)(0)$ & $\C h_2(-\frac{1}{3})\otimes v$ & $\frac{2}{45}$ &
$-\frac{4}{81}\sqrt{-3}$ \\
$W_T(\tau)(1)$ & $\C h_1(-\frac{2}{3})\otimes v$ &
$\frac{2}{45}+\frac{1}{3}$ &
$-\frac{22}{81}\sqrt{-3}$\\
$W_T(\tau)(2)$ & $\C h_2(-\frac{1}{3})^3\otimes v$ &
$\frac{2}{45}+\frac{2}{3}$ &
$\frac{176}{81}\sqrt{-3}$
\end{tabular}
\end{center}
\end{table}

\subsection{Irreducible $M^\tau$-modules in $\tau^2$-twisted
$M$-modules}

Finally, we find $6$ irreducible $M^\tau$-modules in
$\tau^2$-twisted $M$-modules. The argument is parallel to that in
Subsection 4.2. Instead of $\tau$, we take $\tau^2$. Thus we
follow \cite{DL} with $\nu=\tau^2$. Set $h_1'=h_2$, $h_2' =h_1$,
and
\begin{equation*}
\h_{(n)}' = \{\al \in \h\, |\,\tau^2\al=\xi^n\al\}.
\end{equation*}
Then $\h_{(0)}' =0$, $\h_{(1)}' =\C h_1'$, and $\h_{(2)}' =\C
h_2'$. Consider a split central extension
\begin{equation*}
1 \longrightarrow \la\kappa_6\ra \longrightarrow \hat{L}_{\tau^2}
\overset{-}{\longrightarrow} L \longrightarrow 1
\end{equation*}
and choose linear characters $\chi'_j : \hat{L}_{\tau^2}/K \to
\C^\times$, $j=0,1,2$, such that
\begin{equation*}
\chi'_j(\kappa_6) =\xi_6, \quad \chi'_j(\kappa_3
e^{\beta_1}K)=\xi^j, \quad \text{and}\quad \chi'_j(\kappa_3
e^{-\beta_1}K)=\xi^{-j},
\end{equation*}
where $K=\{ a^{-1} \tau^2(a)\,|\,a\in \hat{L}_{\tau^2}\}$. Let
$T_{\chi_j'}$ be the one dimensional $\hat{L}_{\tau^2}/K$-module
affording the character $\chi_j'$. Then the irreducible
$\tau^2$-twisted $V_L$-module associated with $T_{\chi_j'}$ is
\begin{equation*}
V_L^{T_{\chi_j'}}(\tau^2)=S[\tau^2]\otimes T_{\chi_j'}.
\end{equation*}
As a vector space $S[\tau^2]$ is isomorphic to a polynomial
algebra with variables $h_1'(1/3+n)$ and $h_2'(2/3+n)$, $n \in
\Z$. The weight on $S[\tau^2]$ is given by $\wt 1 = 1/9$ and $\wt
h_j'(j/3+n)=-j/3-n$. Moreover, $\wt v=0$ for $v \in T_{\chi_j'}$.
Set
\begin{align*}
M_T(\tau^2)&=\{ u\in V_L^{T_{\chi_0'}}(\tau^2)\,|\,(\om^3)_1 u=0
\},\\
W_T(\tau^2)&=\{ u\in V_L^{T_{\chi_0'}}(\tau^2)\,|\,(\om^3)_1
u=\frac{2}5 u \}.
\end{align*}
Then $M_T(\tau^2)$ and $W_T(\tau^2)$ are the inequivalent
irreducible $\tau^2$-twisted $M$-modules. Furthermore, we have
\begin{align*}
V_L^{T_{\chi_0'}}(\tau^2) & \cong M_T(\tau^2)\otimes
\left(L(\frac{4}{5},0)+L(\frac{4}{5},3)\right) \oplus
W_T(\tau^2)\otimes
\left(L(\frac{4}{5},\frac{2}{5})+L(\frac{4}{5},\frac{7}{5})\right),\\
V_L^{T_{\chi_1'}}(\tau^2) & \cong V_L^{T_{\chi_2'}}(\tau^2) \cong
M_T(\tau^2)\otimes L(\frac{4}{5},\frac{2}{3}) \oplus
W_T(\tau^2)\otimes L(\frac{4}{5},\frac{1}{15})
\end{align*}
as $\tau^2$-twisted $M\otimes \Vir(\om^3)$-modules. The character
of $M_T(\tau^2)$ or $W_T(\tau^2)$ is equal to that of $M_T(\tau)$
or $W_T(\tau)$ respectively. For $\epsilon=0,1,2$, let
\begin{align*}
M_T(\tau^2)(\epsilon) &=\{ u \in M_T(\tau^2)\,|\,\tau^2 u
=\xi^\epsilon u\},\\
W_T(\tau^2)(\epsilon) &=\{ u \in W_T(\tau^2)\,|\,\tau^2
u=\xi^\epsilon u\}.
\end{align*}
Those $6$ modules for $M^\tau$ are inequivalent irreducible
modules by \cite[Theorem 2]{MT}. Their top levels and the
eigenvalues for the action of $L^{\tau^2}(0)=\om_1$ and
$J^{\tau^2}(0)=J_2$ are collected in Table \ref{toptau2}.

\begin{table}[h]
\begin{center}
\caption{irreducible $M^\tau$-modules in $M_T(\tau^2)$ and
$W_T(\tau^2)$} \label{toptau2}
\begin{tabular}{l|l|l|l}
irred. module & top level & $L^{\tau^2}(0)$ & $J^{\tau^2}(0)$\\
\hline $M_T(\tau^2)(0)$ & $\C 1\otimes v$ & $\frac{1}{9}$ &
$-\frac{14}{81}\sqrt{-3}$\\
$M_T(\tau^2)(1)$ & $\C h_2'(-\frac{1}{3})^2\otimes v$ &
$\frac{1}{9}+\frac{2}{3}$ &
$\frac{238}{81}\sqrt{-3}$\\
$M_T(\tau^2)(2)$ & $\C (\frac{4}{3}h_1'(-\frac{2}{3})^2\otimes v +
h_2'(-\frac{1}{3})^4\otimes v)$ & $\frac{1}{9}+\frac{4}{3}$ &
$-\frac{374}{81}\sqrt{-3}$\\
$W_T(\tau^2)(0)$ & $\C h_2'(-\frac{1}{3})\otimes v$ &
$\frac{2}{45}$ &
$\frac{4}{81}\sqrt{-3}$ \\
$W_T(\tau^2)(1)$ & $\C h_1'(-\frac{2}{3})\otimes v$ &
$\frac{2}{45}+\frac{1}{3}$ &
$\frac{22}{81}\sqrt{-3}$\\
$W_T(\tau^2)(2)$ & $\C h_2'(-\frac{1}{3})^3\otimes v$ &
$\frac{2}{45}+\frac{2}{3}$ & $-\frac{176}{81}\sqrt{-3}$
\end{tabular}
\end{center}
\end{table}

\subsection{Remarks on $20$ irreducible $M^\tau$-modules}

We have obtained $20$ irreducible $M^\tau$-modules in Subsections
4.1, 4.2, and 4.3. Note that the top levels of them are of
dimension one and they can be distinguished by the eigenvalues for
$\om_1$ and $J_2$.

The isometry $\sigma$ of the lattice $(L, \la\cdot,\cdot\ra)$
induces a permutation of order $2$ on those $20$ irreducible
$M^\tau$-modules. Clearly, $\sigma$ leaves $M(0)$ and $W(0)$
invariant and transforms $M_k^c$ (resp. $W_k^c$) into an
irreducible $M^\tau$-module equivalent to $M_k^c$ (resp. $W_k^c$).
Moreover, $\sigma$ interchanges irreducible $M^\tau$-modules as
follows:
\begin{equation}\label{sym2}
\begin{array}{cc}
M(1) \longleftrightarrow M(2), &  W(1) \longleftrightarrow W(2),\\
M_T(\tau)(\epsilon) \longleftrightarrow M_T(\tau^2)(\epsilon),
\quad & \quad W_T(\tau)(\epsilon) \longleftrightarrow
W_T(\tau^2)(\epsilon)
\end{array}
\end{equation}
for $\epsilon=0,1,2$. The top level of $M_T(\tau^2)(\epsilon)$ can
be obtained by replacing $h_j(j/3+n)$ with $h_j'(j/3+n)$ for $j=1,
2$ in the top level of $M_T(\tau)(\epsilon)$. Similar symmetry
holds for $W_T(\tau^2)(\epsilon)$ and $W_T(\tau)(\epsilon)$. The
action of $J(0)$ on the top level of $M_T(\tau^2)(\epsilon)$
(resp. $W_T(\tau^2)(\epsilon)$) is negative of the action on the
top level of $M_T(\tau)(\epsilon)$ (resp. $W_T(\tau)(\epsilon)$).
These symmetries are consequences of the fact that
$\sigma\tau\sigma=\tau^2$ and $\sigma J=-J$.

In \cite{FZ} an infinite series of $2D$ conformal field theory
models with $\Z_3$ symmetry was studied. In the case $p=5$ of
\cite{FZ}, $20$ irreducible representations are discussed
\cite[(5.5)]{FZ}. If we multiply the values $w \begin{pmatrix}n
& m\\ n' & m'\end{pmatrix}$ of \cite[(5.6)]{FZ} by
$\sqrt{-105/2}$, then the pairs
\begin{equation*}
\Big( \Delta \begin{pmatrix}n & m\\ n' & m'\end{pmatrix},
\sqrt{-105/2}\, w \begin{pmatrix}n & m\\ n' & m'\end{pmatrix}
\Big)
\end{equation*}
coincide with the pairs of the eigenvalues for $\om_1$ and $J_2$
of the top levels of the $20$ irreducible $M^\tau$-modules listed
in Tables \ref{top1}, \ref{toptau}, and \ref{toptau2}. Here
$\Delta
\begin{pmatrix}n & m\\ n' & m'\end{pmatrix}$ is given by
\cite[(1.3)]{FZ}.

\section{Classification of irreducible modules for $M^\tau$}

We show in this section that the $20$ irreducible modules discussed in
Section 4 are all the inequivalent irreducible modules for
$M^\tau$. This is achieved
by determining the Zhu algebra
$A(\W)$ of the vertex operator subalgebra $\W$ in $M^\tau$
generated by $\om$ and $J$. It turns out that $A(\W)$ is
isomorphic to a quotient algebra of the polynomial algebra $\C[x,
y]$ with two variables $x$ and $y$ by a certain ideal $\I$ and
that $A(\W)$ is of dimension $20$. We will also prove that
$M^\tau = \W$ and $\W$ is rational.

As in Theorem 3.3, let $L(n)=\om_{n+1}$ and $J(n)=J_{n+2}$ for $n
\in \Z$. The action of those operators on the vacuum vector $\1$
is such that
\begin{align}
L(n)\1&=0 \quad \mbox{for} \quad n \ge -1,& L(-2)\1 &= \om,
\label{Lv}\\
J(n)\1&=0 \quad \mbox{for} \quad n \ge -2,& \qquad J(-3)\1 &= J,
\label{Jv}
\end{align}

\subsection{A spanning set for $\W$}
For a vector expressed in the form $u^1_{n_1}\cdots u^k_{n_k}\1$
with $u^i \in \{ \om, J\}$ and $n_i \in \Z$, we denote by $l_\om
(u^1_{n_1}\cdots u^k_{n_k}\1)$ or $l_J (u^1_{n_1}\cdots
u^k_{n_k}\1)$ the number of $i$, $1 \le i \le k$ such that $u^i =
\om$ or $u^i=J$ respectively. We shall call these numbers the
$\om$-length or the $J$-length of the expression $u^1_{n_1}\cdots
u^k_{n_k}\1$. Since each vector in $\W$ is not necessarily
expressed uniquely in such a form, the $\om$-length and the
$J$-length are not defined for a vector. They depend on a specific
expression in the form $u^1_{n_1}\cdots u^k_{n_k}\1$.

\begin{lem}\label{l5.1} Let the $\om$-length and the $J$-length of
$u^1_{n_1}\cdots
u^k_{n_k}\1$ be $s$ and $t$ respectively. Then $u^1_{n_1}\cdots
u^k_{n_k}\1$ can be written as a linear combination of vectors of
the form
\begin{equation*}
L(-m_1)\cdots L(-m_p)J(-n_1)\cdots J(-n_q)\1
\end{equation*}
such that

\medskip
$(1)$ $m_1 \ge \cdots \ge m_p \ge 2, \quad n_1 \ge \cdots \ge n_q
\ge 3$,

\medskip
$(2)$ $q \le t$,

\medskip
$(3)$ $p+q \le s+t$,

\medskip
$(4)$ $m_1 + \cdots + m_p + n_1 + \cdots + n_q = \wt
(u^1_{n_1}\cdots u^k_{n_k}\1)$.
\end{lem}

\prf We proceed by induction on $t$. If $t=0$, the assertion
follows from the commutation relation \eqref{LL} and the action of
$L(n)$ on the vacuum vector \eqref{Lv}.

Suppose the assertion holds for the case where the $J$-length of
$u^1_{n_1}\cdots u^k_{n_k}\1$ is at most $t-1$ and consider the
case where the $J$-length is $t$. By \eqref{LJ}, we can replace
$J(-n)L(-m)$ with $L(-m)J(-n)$ or $J(-m-n)$. Hence we may assume
that $u^1_{n_1}\cdots u^k_{n_k}\1$ is of the form
\begin{equation}\label{LJtemp}
L(-m_1)\cdots L(-m_s)J(-n_1)\cdots J(-n_t)\1
\end{equation}
for some $m_i, n_j \in \Z$.

By \eqref{Jv}, we may assume that $n_t \ge 3$. Suppose $n_i <
n_{i+1}$ for some $i$. Then by the commutation relation
\eqref{JJ}, the vector \eqref{LJtemp} can be written as a linear
combination of the vectors which are obtained by replacing
$J(-n_i)J(-n_{i+1})$ with (i) $J(-n_{i+1})J(-n_i)$, (ii)
$L(-n_i-n_{i+1})$, (iii) $L(k)L(-n_i-n_{i+1}-k)$ or
$L(-n_i-n_{i+1}-k)L(k)$ for some $k \in \Z$, or (iv) a constant.
In Cases (ii), (iii), or (iv), we get an expression whose
$J$-length is at most $t-2$, and so we can apply the induction
hypothesis. Therefore, in \eqref{LJtemp} we may assume that $n_1
\ge \cdots \ge n_t \ge 3$.

Now we argue by induction on the $\om$-length $s$ of the
expression \eqref{LJtemp}. If $s=0$, the assertion holds. Suppose
the assertion holds for the case where the $\om$-length is at most
$s-1$. By \eqref{LJ}, we can replace $L(-m_s)J(-n_1)$ with (i)
$J(-n_1)L(-m_s)$ or (ii) $J(-m_s - n_1)$. In Csae (ii), we get an
expression of $\om$-length at most $s-1$, so that we can apply the
induction hypothesis. Arguing similarly, we can reach
\begin{equation*}
L(-m_1)\cdots L(-m_{s-1})J(-n_1)\cdots J(-n_t)L(-m_s)\1.
\end{equation*}

Hence we may assume that $m_s \ge 2$ by \eqref{Lv}. Suppose $m_i <
m_{i+1}$ for some $i$. Then by \eqref{LL}, the vector
\eqref{LJtemp} can be written as a linear combination of the
vectors which are obtained by replacing $L(-m_i)L(-m_{i+1})$ with
(i) $L(-m_{i+1})L(-m_i)$, (ii) $L(-m_i - m_{i+1})$, or (iii) a
constant. Since Cases (ii) or (iii) yields an expression whose
$\om$-length is at most $s-1$, we can apply the induction
hypothesis. This completes the proof. \quad \qed

\bigskip
A vector of the form
\begin{equation}\label{nf}
L(-m_1)\cdots L(-m_p)J(-n_1)\cdots J(-n_q)\1
\end{equation}
with $m_1 \ge \cdots \ge m_p \ge 2$ and $n_1 \ge \cdots \ge n_q
\ge 3$ will be called of normal form.

\begin{cor}\label{coro1} $\W$ is spanned by the vectors of normal form
\begin{equation*}
L(-m_1)\cdots L(-m_p)J(-n_1)\cdots J(-n_q)\1
\end{equation*}
with $m_1 \ge \cdots \ge m_p \ge 2$, $n_1 \ge \cdots \ge n_q \ge
3$, $p=0,1,2,\ldots$, and $q=0,1,2,\ldots$.
\end{cor}

\prf As a vector space $\W$ is spanned by the vectors
$u^1_{n_1}\cdots u^k_{n_k}\1$ with $u^i \in \{ \om, J\}$, $n_i \in
\Z$, and $k=0,1,2,\ldots$. Hence the assertion follows from Lemma
5.1. \quad \qed

\bigskip
A similar argument for a spanning set can be found in
\cite[Section3]{DN1}. See also \cite[Section 2.2]{BMP}.

\begin{rmk}\label{r1} Let $U$ be an admissible $\W$-module generated
by $u\in U$ such that $L(n)u=J(n)u=0$ for $n>0$ and $L(0)u=hu$,
$J(0)u=ku$ for some $h, k \in \C$. It can be proved in a same way
that $U$ is spanned by
\begin{equation*}
L(-m_1)\cdots L(-m_p)J(-n_1)\cdots J(-n_q)u
\end{equation*}
with $m_1 \ge \cdots \ge m_p \ge 1$, $n_1 \ge \cdots \ge n_q \ge
1$, $p=0,1,2,\ldots$, and $q=0,1,2,\ldots$.
\end{rmk}

\subsection{A singular vector $\bv^{12}$}

A singular vector $v$ of weight $h$ for $\W$ is by definition a
vector $v$ which satisfies

\medskip
$(1)$ $L(0)v=hv$,

\medskip
$(2)$ $L(n)v = 0$ and $J(n)v=0$ for $n \ge 1$.

\medskip\noindent Note that $v$ is not necessarily an eigenvector
for $J(0)$. By commutation relations \eqref{LL} and \eqref{LJ}, it
is easy to show that the condition $(2)$ holds if $v$ satisfies

\medskip
$(2)'$ $L(1)v=L(2)v=J(1)v=0$.

\medskip
We consider $\W$ as a space spanned by the vectors of the form
\eqref{nf} . The weight of such a vector is
$m_1+\cdots+m_p+n_1+\cdots+n_q$. Let $v$ be a linear combination
of the vectors of the form \eqref{nf} of weight $h$. For example,
there are $76$ vectors of the form \eqref{nf} of weight $12$. We
use the conditions \eqref{Lv} and \eqref{Jv} and the commutation
relations \eqref{LL}, \eqref{LJ}, and \eqref{JJ} to compute
$L(1)v$, $L(2)v$, and $J(1)v$. This computation was done by a
computer algebra system Risa/Asir. The result is as follows.

\begin{lem} Let $v$ be a linear combination of the vectors of the
form \eqref{nf} of weight $h$. Under the conditions \eqref{Lv} and
\eqref{Jv} and the commutation relations \eqref{LL}, \eqref{LJ},
and \eqref{JJ}, we have $L(1)v=L(2)v=J(1)v=0$ only if $v=0$ in the
case $h \le 11$. In the case $h=12$, there exists a unique, up to
scalar multiple, linear combination $\bv^{12}$ which satisfies
$L(1)\bv^{12}=L(2)\bv^{12}=J(1)\bv^{12}=0$. The explicit form of
$\bv^{12}$ is given in Appendix A. We also have $J(0)\bv^{12} =
0$.
\end{lem}

We only use the conditions \eqref{Lv} and \eqref{Jv} and the
commutation relations \eqref{LL}, \eqref{LJ}, and \eqref{JJ} to
obtain $\bv^{12}$ in the above computation. Since we consider $\W$
inside the lattice vertex operator algebra $V_L$, there might
exist some nontrivial relations among the vectors of the form
\eqref{nf} which are not known so far. This ambiguity will be
removed in Subsection 5.3.

\subsection{An invariant positive definite hermitian form on
$V_L$}

It is well known that the vertex operator algebra constructed from
any positive definite even lattice as in \cite{FLM} possesses a
positive definite hermitian form which is invariant in a certain
sense (\cite{FHL, FLM, Li, M3}). Following \cite[Section 2.5]{M3},
we review it for our $V_L$.

Set $\widet{L}(n) = \widet{\om}_{n+1}, n \in \Z$, where
$\widet{\om}$ is the Virasoro element of $V_L$. Then
$\widet{L}(1)(V_L)_{(1)} = 0$ and $(V_L)_{(0)}$ is one
dimensional. Thus by \cite[Theorem 3.1]{Li}, there is a unique
symmetric invariant bilinear form $(\,\cdot\, , \,\cdot\,)$ on
$V_L$ such that $(\1,\1)=1$. That the form is invariant means
\begin{equation}\label{inv}
(Y(u,z)v,w)=(v,Y(e^{z\widet{L}(1)}(-z^{-2})^{\widet{L}(0)}u,
z^{-1})w)
\end{equation}
for $u,v,w \in V_L$. The value $(u,v)$ is determined by
\begin{align}
(\1,\1) &= 1, \label{form1}\\
(u,v) &= \Res_z
z^{-1}(\1,Y(e^{z\widet{L}(1)}(-z^{-2})^{\widet{L}(0)}u,z^{-1})v).
\label{form2}
\end{align}

>From \eqref{inv}, we see that
$(\widet{L}(n)u,v)=(u,\widet{L}(-n)v)$. In case of $n=0$, this
implies $((V_L)_{(m)}, (V_L)_{(n)}) = 0$ if $m \ne n$. For $\al
\in L$ and $u,v \in V_L$,
\begin{equation}\label{form3}
\begin{split}
(\al(n)u,v) &= \Res_z z^n (Y(\al(-1),z)u,v)\\
&= -(u,\al(-n)v). \end{split}
\end{equation}

Furthermore, for $\al,\beta \in L$ we have
\begin{equation}\label{form4}
(e^\al,e^\beta) = \delta_{\alpha+\beta, 0}.
\end{equation}

Note that $(-1)^{\la\al,\al\ra/2} = 1$ since $\al \in L$. Consider
an $\R$-form $V_{L,\R}$ of $V_L$ as in \cite[Section 12.4]{FLM}.
That is, let $M(1)_{\R} = \R [ \al(n);\,\al \in L,\,n<0 ]$ and
$V_{L,\R} = M(1)_{\R} \otimes \R [L]$. Then $\C \otimes_{\R} V_{L,
\R} = V_L$. Moreover, $V_{L,\R}$ is invariant under the
automorphism $\theta$. Let $V_{L,\R}^{\pm} = \{ v \in V_{L,\R}
\,|\, \theta v=\pm v\}$. We shall show that the form $(\,\cdot\, ,
\,\cdot\,)$ is positive definite on $V_{L,\R}^+$ and negative
definite on $V_{L,\R}^-$. Indeed, let $\{ \gamma_1,\gamma_2\}$ be
an orthonormal basis of $\R \otimes_{\Z} L$. Then using
\eqref{form3} and \eqref{form4} we can verify that
\begin{equation*}
(\gamma_{i_1}(-m_1)\cdots \gamma_{i_p}(-m_p)e^{\al},
\gamma_{j_1}(-n_1)\cdots \gamma_{j_q}(-n_q)e^{\beta}) \ne 0
\end{equation*}
only if $\gamma_{i_1}(-m_1)\cdots \gamma_{i_p}(-m_p) =
\gamma_{j_1}(-n_1)\cdots \gamma_{j_q}(-n_q)$ in $M(1)_{\R}$ and
$\al+\beta =0$. Furthermore,
\begin{equation}\label{form5}
(\gamma_{i_1}(-m_1)\cdots \gamma_{i_p}(-m_p)e^{\al},
\gamma_{i_1}(-m_1)\cdots \gamma_{i_p}(-m_p)e^{-\al}) = (-1)^p
\cdot (\mbox{a positive integer}).
\end{equation}

We can choose a basis of $V_{L,\R}^+$ consisting of vectors of the
form
\begin{align*}
\gamma_{i_1}(-m_1)\cdots \gamma_{i_p}(-m_p)(e^{\al}+e^{-\al}), &
\quad p \mbox{\ even},\quad \al \in L,\\
\gamma_{i_1}(-m_1)\cdots \gamma_{i_p}(-m_p)(e^{\al}-e^{-\al}), &
\quad p \mbox{\ odd},\quad 0 \ne \al \in L.
\end{align*}

By \eqref{form4} and \eqref{form5}, these vectors are mutually
orthogonal and the square length of each of them is a positive
integer. Hence the form $(\,\cdot\, , \,\cdot\,)$ is positive
definite on $V_{L,\R}^+$. Likewise, we see that the form
$(\,\cdot\, , \,\cdot\,)$ is negative definite on $V_{L,\R}^-$.

We also have $(V_{L,\R}^+, V_{L,\R}^-)=0$. Thus the form
$(\,\cdot\, , \,\cdot\,)$ is positive definite on $V_{L,\R}^+ +
\sqrt{-1}V_{L,\R}^-$. The $\R$-vector space $V_{L,\R}^+ +
\sqrt{-1}V_{L,\R}^-$ is an $\R$-form of $V_L$ since $V_L=\C
\otimes_{\R} (V_{L,\R}^+ + \sqrt{-1}V_{L,\R}^-)$. Note that it is
invariant under the component operators $u_n$ of $Y(u,z)$ for $u
\in V_{L,\R}^+$.

Define a hermitian form $((\,\cdot\, , \,\cdot\,))$ on $V_L$ by
$((\lambda u,\mu v)) = \lambda\overline{\mu}(u,v)$ for $\lambda,
\mu \in \C$ and $u,v \in V_{L,\R}^+ + \sqrt{-1} V_{L,\R}^-$. Then
the hermitian form $((\,\cdot\, , \,\cdot\,))$ is positive
definite on $V_L$ and invariant under $V_{L,\R}^+$, that is,

\begin{equation}\label{form6} ((Y(u,z)v,w)) =
((v,Y(e^{z\widet{L}(1)}(-z^{-2})^{\widet{L}(0)}u,z^{-1})w))
\end{equation}
for $u \in V_{L,\R}^+$ and $v,w \in V_L$.

Using the hermitian form $((\,\cdot\, , \,\cdot\,))$, we can show
that $V_L$ is semisimple as a $\W$-module and that $\W$ is a
simple vertex operator algebra. Note that $\widet{L}(n)v=L(n)v$
for $v \in M$. Note also that $V_{L,\R}^+$ contains $\om$ and $J$.
Then by \eqref{form6},
\begin{align}
((L(n)u, v)) &= ((u, L(-n)v)), \label{form7}\\
((J(n)u, v)) &= -((u, J(-n)v)) \label{form8}
\end{align}
for $n \in \Z$ and $u,v \in V_L$.

Let $U$ be a $\W$-submodule. Denote by $U^{\perp}$ the orthogonal
complement of $U$ in $V_L$ with respect to $((\,\cdot\,
,\,\cdot\,))$. Then $V_L = U \oplus U^{\perp}$ since $((\,\cdot\,
, \,\cdot\,))$ is positive definite. Moreover, $U^{\perp}$ is also
a $\W$-submodule by \eqref{form7} and \eqref{form8}. Thus we
conclude that

\begin{thm}
$V_L$ is semisimple as a $\W$-module.
\end{thm}

Since the weight $0$ subspace $\C\1$ of $\W$ is one dimensional
and since $\W$ is generated by $\1$ as a $\W$-module, we have

\begin{thm}
$\W$ is a simple vertex operator algebra.
\end{thm}

Then there is no singular vector in $\W$ of positive weight.
Hence,

\begin{cor}
The singular vector $\bv^{12}=0$.
\end{cor}

\subsection{The Zhu algebra $A(\W)$}

Based on the properties of $\W$ we have obtained so far, we shall
determine the Zhu algebra $A(\W)$ of $\W$. First we review some
notations and formulas for the Zhu algebra $A(V)$ of an arbitrary
vertex operator algebra $(V,Y,\1,\om)$. The standard reference is
\cite[Section 2]{Z}.

For $u,v \in V$ with $u$ being homogeneous, define two binary
operations
\begin{align}
u \ast v &= \Res_z\Big(\frac{(1+z)^{\wt u}}{z}Y(u,z)v\Big) =
\sum_{i=0}^{\infty} \binom{\wt u}{i}u_{i-1}v, \label{prd1}\\
u \circ v &= \Res_z\Big(\frac{(1+z)^{\wt u}}{z^2}Y(u,z)v\Big) =
\sum_{i=0}^{\infty} \binom{\wt u}{i}u_{i-2}v. \label{prd2}
\end{align}

We extend $\ast$ and $\circ$ for arbitrary $u, v \in V$ by
linearity. Let $O(V)$ be the subspace of $V$ spanned by all $u
\circ v$ for $u,v \in V$. By a theorem of Zhu \cite{Z}, $O(V)$ is
a two-sided ideal with respect to the operation $\ast$. Thus it
induces an operation on $A(V)=V/O(V)$. Denote by $[v]$ the image
of $v \in V$ in $A(V)$. Then $[u]\ast [v] = [u\ast v]$ and $A(V)$
is an associative algebra by this operation. Moreover, $[\1]$ is
the identity and $[\om]$ is in the center of $A(V)$. We denote by
$[u]^{\ast p}$ the product of $p$ copies of $[u]$ in $A(V)$. For
$u,v \in V$, we write $u \sim v$ if $[u] = [v]$. For $f,g \in \End
V$, we write $f \sim g$ if $fv \sim gv$ for all $v \in V$. We need
some formulas from \cite{Z}.
\begin{equation}\label{gen}
\Res_z\Big(\frac{(1+z)^{\wt (u) + m}}{z^{2+n}}Y(u,z)v\Big) =
\sum_{i=0}^{\infty} \binom{\wt (u) + m}{i}u_{i-n-2}v \in O(V)
\end{equation}
for $n \ge m \ge 0$ and
\begin{equation}\label{vu}
v \ast u \sim \Res_z\Big(\frac{(1+z)^{\wt (u) - 1}}{z}Y(u,z)v\Big)
= \sum_{i=0}^{\infty} \binom{\wt (u) - 1}{i}u_{i-1}v.
\end{equation}

Moreover (see \cite{Wa1}),
\begin{equation}\label{virn}
L(-n) \sim (-1)^n \big\{ (n-1)\big( L(-2)+L(-1)\big) + L(0)\big\}
\end{equation}
for $n \ge 1$ and
\begin{equation}\label{virp}
[\om] \ast [v] = [(L(-2)+L(-1))v].
\end{equation}

It follows from \eqref{virn} and \eqref{virp} that
\begin{equation}\label{f1}
[L(-n)u] = (-1)^n (n-1) [\om]\ast [u] + (-1)^n [L(0)u]
\end{equation}
for $n \ge 1$.

For a homogeneous $u \in V$, set $o(u)=u_{\wt (u) - 1}$, which is
the weight zero component operator of $Y(u,z)$. Extend $o(u)$ for
arbitrary $u \in V$ by linearity. We call a module in the sense of
\cite{Z} as an admissible module as in \cite{DLM}. If $U =
\oplus_{n=0}^{\infty} U(n)$ is an admissible $V$-module with
$U(0)\ne 0,$ then
$o(u)$ acts on its top level $U(0)$. Zhu's theory \cite{Z} says

\medskip
(1) $o(u)o(v) = o(u \ast v)$ as operators on the top level $U(0)$
and $o(u)$ acts as $0$ on $U(0)$ if $u \in O(V)$. Thus $U(0)$ is an
$A(V)$-module, where $[u]$ acts on $U(0)$ as $o(u)$.

\medskip
(2) The map $U \mapsto U(0)$ is a bijection between the set of
equivalence classes of irreducible admissible $V$-modules and the
set of equivalence classes of irreducible $A(V)$-modules.

\medskip
We now return to $\W$. Since $\wt J = 3$, we have
\begin{equation}\label{f2}
[J(-n-4)v] = -3[J(-n-3)v] - 3[J(-n-2)v] - [J(-n-1)v]
\end{equation}
for $v \in \W$ and $n \ge 0$ by \eqref{gen}.

\begin{lem}
The image $[L(-m_1)\cdots L(-m_p)J(-n_1)\cdots J(-n_q)\1]$ of the
vector \eqref{nf} with $m_1 \ge \cdots \ge m_p \ge 2$ and $n_1 \ge
\cdots \ge n_q \ge 3$ in $A(\W)$ is contained in
\begin{equation*}
\spn \{ [\om]^{\ast s} \ast [J]^{\ast t} \,|\,0 \le s,\,0\le t \le
q,\,2s+3t \le m_1+\cdots+m_p+n_1+\cdots+n_q \}.
\end{equation*}

In particular, $A(\W)$ is commutative and every element of $A(\W)$
is a polynomial in $[\om]$ and $[J]$.
\end{lem}

\prf We proceed by induction on the $J$-length $q$. By a repeated
use of \eqref{f1}, we see that $[L(-m_1)\cdots
L(-m_p)J(-n_1)\cdots J(-n_q)\1]$ is a linear combination of
$[\om]^{\ast s} \ast [J(-n_1)\cdots J(-n_q)\1]$, $0 \le s \le p$.
Thus the assertion holds if $q=0$.

Suppose the assertion holds for vectors of normal form with
$J$-length at most $q-1$ and consider $[J(-n_1)\cdots J(-n_q)\1]$.
Let $v=J(-n_1)\cdots J(-n_q)\1$ and $u=J(-n_2)\cdots J(-n_q)\1$,
so that $v=J(-n_1)u$. we proceed by induction on the weight. The
vector of the smallest weight is the case $n_1=3$. In this case
$v=J(-3)^q\1$ and $u=J(-3)^{q-1}\1$. Since $v=J_{-1}u$, it follows
from \eqref{prd1} that
\begin{equation*}
[v]=[J]\ast [u] - 3[J(-2)u] - 3[J(-1)u] - [J(0)u].
\end{equation*}

The weight of $J(-n)u$, $0 \le n \le 2$, is less than $\wt v$. By
Lemma 5.1, each of these three vectors is a linear combination of
vectors of normal form with $J$-length at most $q-1$. Then we can
apply the induction hypothesis on $J$-length and the assertion
holds if $n_1=3$. Assume that $n_1 \ge 4$. By \eqref{f2},
$[v]=[J(-n_1)u]$ is a linear combination of $[J(-n)u]$, $n_1-3 \le
n \le n_1-1$. The weight of $J(-n)u$, $n_1-3 \le n \le n_1-1$, is
less than $\wt v$. Hence by Lemma 5.1, these three vectors are
linear combinations of vectors of normal form with $J$-length at
most $q$ and weight less than $\wt v$. The induction is complete.
\quad \qed

\bigskip
The image $[L(-m_1)\cdots L(-m_p)J(-n_1)\cdots J(-n_q)\1]$ of the
vector of normal form \eqref{nf} with $m_1 \ge \cdots \ge m_p \ge
2$ and $n_1 \ge \cdots \ge n_q \ge 3$ in $A(\W)$ can be written
explicitly as a polynomial in $[\om]$ and $[J]$ by the following
algorithm.

Since $A(\W)$ is commutative, it follows from \eqref{vu} that
\begin{equation}\label{f3}
[J(-3)v] = [J]\ast [v] - 2[J(-2)v] - [J(-1)v]
\end{equation}
for $v \in \W$. Now we use \eqref{f1}, \eqref{f2}, and \eqref{f3}.
Although $J(-n-4)v$ is of normal form, the vectors $J(-n-3)v$,
$J(-n-2)v$, and $J(-n-1)v$ in \eqref{f2} may not be of normal
form. However, the weight of any of these three vectors is less
than the weight of $J(-n-4)v$, and so we can apply the argument in
the proof of Lemma 5.1. A similar discussion is also needed for
the formula \eqref{f3}. Thus the algorithm is by induction on the
weight. We use formulas \eqref{f1}, \eqref{f2}, \eqref{f3} and
apply Lemma 5.1, that is, use the commutation relations
\eqref{LL}, \eqref{LJ}, \eqref{JJ} and the conditions \eqref{Lv}
and \eqref{Jv}. By induction on the weight and a repeated use of
those formulas and conditions, we can write explicitly the image
of the vector \eqref{nf} in $A(\W)$ as a polynomial in $[\om]$ and
$[J]$.

Consider the algebra homomorphism
\begin{equation*}
\C [x,y] \longrightarrow A(\W); \quad x \longmapsto [\om], \quad y
\longmapsto [J]
\end{equation*}
of the polynomial algebra $\C [x,y]$ with two variables $x,y$ onto
$A(\W)$. Denote its kernel by $\I$. Then $\C[x,y]/\I \cong A(\W)$.
We shall consider $\bv^{12}$, $J(-1)\bv^{12}$, $J(-2)\bv^{12}$,
and $J(-1)^2\bv^{12}$. These vectors are described explicitly as
linear combinations of vectors of normal form in Appendix A. Their
images $[\bv^{12}]$, $[J(-1)\bv^{12}]$, $[J(-2)\bv^{12}]$, and
$[J(-1)^2\bv^{12}]$ can be written as polynomials in $[\om]$ and
$[J]$ by the above mentioned algorithm. The results are given in
Appendix B. Let $F_i(x,y) \in \C[x,y]$, $1 \le i \le 4$, be the
polynomials which are obtained by replacing $[\om]$ with $x$ and
$[J]$ with $y$ in the polynomials given in Appendix B. Since
$\bv^{12}=0$ by Corollary 5.6, $F_i(x,y)$'s are contained in $\I$.
Let $\I'$ be the ideal in $\C[x,y]$ generated by $F_i(x,y), 1 \le
i \le 4$.

The primary decomposition of $\I'$ is $\I' = \cap_{i=1}^{20}
\CP_i$, where $\CP_i$, $1 \le i \le 20$ are
\begin{equation}\label{ideal}
\begin{array}{ll}
\la x,\, y \ra, &
\la 5x-8,\, y\ra,\\
\la 2x-1,\, y\ra, &
\la 10x-1,\, y\ra,\\
\la x-2,\, y-12\sqrt{-3}\ra, &
\la x-2,\, y+12\sqrt{-3}\ra,\\
\la 5x-3,\, y-2\sqrt{-3}\ra,&
\la 5x-3,\, y+2\sqrt{-3}\ra,\\
\la 9x-1,\, 81y-14\sqrt{-3}\ra, &
\la 9x-1,\,81y+14\sqrt{-3}\ra,\\
\la 9x-7,\, 81y-238\sqrt{-3}\ra,&
\la 9x-7,\, 81y+238\sqrt{-3}\ra,\\
\la 9x-13,\, 81y-374\sqrt{-3}\ra, & \la
9x-13,\,81y+374\sqrt{-3}\ra,\\
\la 45x-2,\, 81y-4\sqrt{-3}\ra, &
\la 45x-2,\, 81y+4\sqrt{-3}\ra,\\
\la 45x-17,\, 81y-22\sqrt{-3}\ra, & \la
45x-17,\,81y+22\sqrt{-3}\ra,\\
\la 45x-32,\,81y-176\sqrt{-3}\ra,& \la
45x-32,\,81y+176\sqrt{-3}\ra,
\end{array}
\end{equation}

These primary ideals correspond to the $20$ irreducible
$M^{\tau}$-modules listed in Tables 1, 3, and 4 in Section 4. The
correspondence is given by substituting $x$ and $y$ with the
eigenvalues for $L(0)$ and $J(0)$ on the top levels of $20$
irreducible modules. The eigenvalues are the zeros of those
primary ideals. Note that the $20$ pairs of those eigenvalues for
$L(0)$ and $J(0)$ on the top levels are different from each other.
Since the top levels of the $20$ irreducible $M^{\tau}$-modules
are one dimensional and since $\W$ is contained in $M^{\tau}$,
there are at least $20$ inequivalent irreducible $\W$-modules
whose top levels are the same as those of irreducible
$M^{\tau}$-modules. Hence by Zhu's theory \cite{Z}, we conclude
that $\I=\I'$ and $A(\W) \cong \oplus_{i=1}^{20} \C[x,y]/\CP_i$.
In particular, $\W$ has exactly $20$ inequivalent irreducible
modules.

If $\W \ne M^{\tau}$, then we can take an irreducible $\W$-module
$U$ in $M^{\tau}$ such that $\W \cap U = 0$ by Theorem 5.4. From
the classification of irreducible $\W$-modules we see that the
smallest weight of $U$ is at most $2$. But we can verify that the
homogeneous subspaces of $\W$ of weight $0$, $1$, and $2$ coincide
with those of $M^{\tau}$. Therefore, $\W = M^{\tau}$.

We have obtained the following theorem.

\begin{thm}\label{tt1}
$(1)$ $M^{\tau} = \W$.

$(2)$ $A(M^{\tau}) \cong \oplus_{i=1}^{20} \C[x,y]/\CP_i$
is a 20-dimensional commutative associative algebra.

$(3)$ There are exactly $20$ inequivalent irreducible
$M^{\tau}$-modules. Their representatives are listed in Tables 1,
3, and 4 in Section 4, namely, $M(\epsilon)$, $W(\epsilon)$,
$M_k^c$, $W_k^c$, $M_T(\tau^i)(\epsilon)$, and
$W_T(\tau^i)(\epsilon)$ for $\epsilon = 0,1,2$ and $i=1,2$.
\end{thm}

\begin{rmk}
The explicit description of $\bv^{12}$, $J(-1)\bv^{12}$,
$J(-2)\bv^{12}$, and $J(-1)^2\bv^{12}$ in Appendix A, the images
of these four vectors in $A(\W)$ in Appendix B, and the primary
ideals \eqref{ideal} were obtained by a computer algebra system
Risa/Asir.
\end{rmk}

\subsection{Rationality of $\W$}

First we prove a result about a general vertex operator algebra.
For this purpose, we need to recall that a vertex operator
algebra $V$ is called $C_2$-cofinite if $V/C_2(V)$
is finite dimensional where $C_2(V)$ is the subspace
of $V$ spanned by $u_{-2}v$ for $u,v\in V.$

\begin{prop}\label{p1}
Let $V=\oplus_{n\geq 0} V_n$ be a $C_2$-cofinite vertex operator
algebra such that $V_0$ is one-dimensional. Assume that $A(V)$ is
semisimple and any $V$-module generated by an irreducible
$A(V)$-module is irreducible. Then $V$ is a rational vertex
operator algebra.
\end{prop}

\prf By the definition of rationality (cf. \cite{dlmnew}), we need
to prove that any admissible $\W$-module $Z$ is completely
reducible. By \cite[Lemma 5.5]{abd}, $Z$ is a direct sum of
generalized eigensapces for $L(0).$ So it is enough to prove that
any submodule generated by a generalized eigenvector for $L(0)$ is
completely reducible. We can assume that $Z$ is generated by a
generalized eigenvector for $L(0).$ Then  $Z=\oplus_{n\geq
0}Z_{\lambda+n}$ for some $\lambda\in \C$ where $Z_{\lambda+n}$ is
the generalized eigenspace for $L(0)$ with eigenvalue $\lambda+n$
and $Z_{\lambda}\ne 0.$ We call $\lambda$ the minimal weight of
$Z.$ By \cite[Theorem 1]{buhl}, each $Z_{\lambda+n}$ is finite
dimensional.

Let $X$ be the submodule of $Z$ generated by $Z_{\lambda}.$ Then
$X$ is completely reducible by the assumption. So we have an exact
sequence
$$0\to X\to Z\to Z/X\to 0$$
of admissible $V$-modules. Let $Z'=\oplus_{n\geq 0}Z_{\lambda+n}^*$
be the graded dual of $Z.$ Then $Z'$ is also an admissible
$V$-module (see \cite{FHL}) and
we have an exact sequence
$$0\to (Z/X)'\to Z'\to X'\to 0$$
of admissible $V$-modules. On the other hand, the $V$-submodule
of $Z'$ generated by $Z_{\lambda}^*$ is isomorphic to $X'.$ As a result
we have $Z'$ is isomorphic to $X'\oplus (Z/X)'.$ This implies
that $Z\cong X\oplus Z/X.$  Clearly, the minimal weight of
$Z/X$ is greater than the minimal weight of $Z.$ Continuing
in this way we prove that $Z$ is a direct sum irreducible modules. \quad
\qed

\medskip
Now we turn our attention to $\W.$

\begin{thm}\label{l5.52} $\W$ is $C_2$-cofinite.
\end{thm}

\prf Note from Corollary \ref{coro1} that $\W$ is spanned by
\begin{equation*}
L(-m_1)\cdots L(-m_p)J(-n_1)\cdots J(-n_q)\1
\end{equation*}
with $m_1 \ge \cdots \ge m_p \ge 2$, $n_1 \ge \cdots \ge n_q \ge
3$, $p=0,1,2,\ldots$, and $q=0,1,2,\ldots$. Then $\W$ is spanned
by $L(-2)^pJ(-3)^q\1$ modulo $C_2(\W).$ It is well known that
$\W/C_2(\W)$ is a commutative associative algebra under the
product $u\cdot v=u_{-1}v$ for $u,v\in \W$ (cf. \cite{Z}). So $\W$
is spanned by $\omega^p\cdot J^q$ modulo $C_2(\W)$ for $p,q\geq
0.$

The key idea to prove that $\W$ is $C_2$-cofinite is to use the
singular vector $\bv^{12}$. By the explicit form of $\bv^{12}$,
$J(-1)\bv^{12}$, and $J(-1)^2\bv^{12}$ in Appendix A, we have the
following relations in $\W/C_2(\W)$.
\begin{eqnarray*}
& &-(59680000/3501)\omega^6 -(184400/1167)\omega^3\cdot J^2+J^4=0,\\
& &-926640\omega^2\cdot J^3-89856000\omega^5 \cdot J=0,\\
& &21565440000\omega^7-680659200\omega^4\cdot J^2-5559840\omega\cdot J^4=0.
\end{eqnarray*}
Multiplying by $\omega^2, J, \omega$ respectively we get
\begin{eqnarray*}
& &-(59680000/3501)\omega^8 -(184400/1167)\omega^5\cdot J^2+\omega^2\cdot J^4=0,\\
& &-926640\omega^2\cdot J^4-89856000\omega^5 \cdot J^2=0,\\
& &21565440000\omega^8-680659200\omega^5\cdot J^2-5559840\omega^2\cdot J^4=0.
\end{eqnarray*}
It follows immediately that
\begin{equation*}
\omega^8=\omega^2\cdot J^4=\omega^5 \cdot J^2=0.
\end{equation*}
Thus
\begin{equation*}
J^8=\left((59680000/3501)\omega^6 +(184400/1167)\omega^3\cdot
J^2\right)^2=0.
\end{equation*}
As a result, $\W/C_2(\W)$ is
spanned by $\omega^p\cdot J^q$ for $0\leq p,q\leq 7,$ as desired.
\quad \qed

\begin{lem}\label{l5.53} Let $U$ be an irreducible $A(\W)$-module. Then any
$\W$-module $Z$ generated by $U$ is irreducible.
\end{lem}

\prf By Theorem \ref{tt1}, $A(\W)$ has exactly twenty irreducible
modules and $\omega$ acts on each irreducible module as a constant
in the set
\begin{equation*}
P=\{0, 2, 8/5,3/5,1/2,1/10, 1/9,1/9+2/3, 1/9+4/3, 2/45,
2/45+1/3,2/45+2/3\}.
\end{equation*}

Let $\omega$ act on $U$ as $\lambda.$ Assume that $\lambda \ne 0,
3/5.$ Then $\lambda$ is maximal in the set $P\cap (\lambda+\Z).$
Let $Z=\oplus_{n\ge0}Z_{\lambda+n}$ and $Z_{\lambda}=U.$ If $Z$ is
not irreducible then $Z$ has a proper submodule $X=\sum_{n\geq
0}X_{\lambda+n_0+n}$ for some $n_0>0$ with $X_{\lambda+n_0}\ne 0$
where $X_{\lambda+m}= X\cap Z_{\lambda+m}.$ So $X_{\lambda+n_0}$
is an $A(\W)$-module on which $\omega$ acts on $\lambda+n_0.$
Since $\lambda+n_0\in P\cap (\lambda+\Z)$ is greater than
$\lambda$ we have a contradiction. This shows that $Z$ must be
irreducible.

It remains to prove the result with $\lambda=0$ or $\lambda=3/5$.
If $\lambda=0,$ then $U\cong \C\1$ and $Z$ is isomorphic to $\W$
(see \cite{Li}). Now let $\lambda=3/5.$ By Theorem \ref{tt1}, $U$
can be either $W(1)_{3/5}$ or $W(2)_{3/5}$ (see Table 1). We can
assume that $U=W(1)_{3/5}$ and the proof for $U=W(2)_{3/5}$ is
similar. In this case $J(0)$ acts on $U$ as $2\sqrt{-3}.$ Let
$U=\C u.$ Then $Z$ is spanned by
\begin{equation*}
L(-m_1)\cdots L(-m_p)J(-n_1)\cdots J(-n_q)u
\end{equation*}
with $m_1 \ge \cdots \ge m_p \ge 1$, $n_1 \ge \cdots \ge n_q \ge
1$, $p=0,1,2,\ldots$, and $q=0,1,2,\ldots$ (see Remark \ref{r1}).
Since $8/5$ is the only number in $P\cap (3/5+\Z)$ greater than
$3/5,$ $Z$ is irreducible if and only if there is no nonzero
vector $v\in Z_{8/5}$ such that $L(1)v=J(1)v=0.$

Note that $Z_{8/5}$ is spanned by $L(-1)u$ and $J(-1)u.$ By
formulas (\ref{LL})-(\ref{JJ}) we see that
\begin{eqnarray*}
& & L(1)L(-1)u=\frac{6}{5}u,\\
& & L(1)J(-1)u=6\sqrt{-3}u,\\
& & J(1)L(-1)u=6\sqrt{-3}u,\\
& & J(1)J(-1)u=\left(\frac{237\times 6}{5}-\frac{48\times 39}{5}\right)u.
\end{eqnarray*}

Now let $v=\alpha L(-1)u+\beta J(-1)u\in Z_{8/5}$ such that
$L(1)v=J(1)v=0$. Then we have a system of linear equations
\begin{eqnarray*}
& &\frac{6}{5}\alpha+6\sqrt{-3}\beta=0,\\
& &6\sqrt{-3}\alpha -90\beta=0.
\end{eqnarray*}
Unfortunately, the system is degenerate and has solutions
$\alpha=-5\sqrt{-3}\beta.$ Thus up to a constant we can assume
that $v=-5\sqrt{-3}L(-1)u+J(-1)u.$

We have to prove that $v=0.$ If $v$ is not zero, then $\C v$ is an
irreducible module for $A(\W)$ on which $\omega$ acts as $8/5.$
Again by Table 1, $J$ must acts on $v$ as $0$. Using (\ref{LJ})
and (\ref{JJ}), we find out that
$J(0)v=-120L(-1)u-8\sqrt{-3}J(-1)u=-8\sqrt{-3}v.$ This implies
that $v=0.$ Clearly we have a contradiction. Thus $Z$ is an
irreducible $\W$-module. \quad \qed

\medskip
Combining Propsotion \ref{p1}, Theorem \ref{l5.52}, and Lemma
\ref{l5.53} together yields
\begin{thm} The vertex operator algebra $\W$ is rational.
\end{thm}

It is proved in \cite{abd} that a rational and $C_2$-cofinite
vertex operator algebra is regular in the sense that any weak module
is a direct sum of irreducible admissible modules. Thus we, in fact,
have proved that $\W$ is also regular.

\section{Characters of irreducible $M^{\tau}$-modules}

We shall describe the characters of the $20$ irreducible
$M^{\tau}$-modules by the characters of irreducible modules for
the Virasoro vertex operator algebras. Throughout this section $z$
denotes a complex number in the upper half plane $\CH$ and $q=\exp
(2\pi\sqrt{-1}z)$. First we recall the character of the
irreducible module $L(c_m,h^{(m)}_{r,s})$ with highest weight
$h^{(m)}_{r,s}$ for the Virasoro vertex operator algebra
$L(c_m,0)$ with central charge $c_m$, where
\begin{align*}
c_m &= 1-\frac{6}{(m+2)(m+3)},\quad m=1,2,\ldots,\\
h^{(m)}_{r,s} &= \frac{\big((m+3)r-(m+2)s\big)^2-1}{4(m+2)(m+3)},
\quad 1\leq s\leq r\leq m+1.
\end{align*}

The character of $L(c_m,h^{(m)}_{r,s})$ is obtained in \cite{R} as
follows:
\begin{equation}\label{eqn:ch}
\ch L(c_m,h^{(m)}_{r,s})=\frac{\sum_{k\in{\mathbb Z}}(q^{b(k)}-
q^{a(k)})}{\prod_{i=1}^{\infty}(1-q^i)},
\end{equation}
where
\begin{align*}
a(k)&=\frac{\big(2(m+2)(m+3)k+(m+3)r+(m+2)s\big)^2-1}{4(m+2)(m+3)},\\
b(k)&=\frac{\big(2(m+2)(m+3)k+(m+3)r-(m+2)s\big)^2-1}
{4(m+2)(m+3)}.
\end{align*}

Define ${\displaystyle \Xi^{(m)}_{r,s}(z) = q^{-c_m/24} \ch
L(c_m,h^{(m)}_{r,s})}$. For $1\leq s\leq r\leq m+1$, the following
transformation formula holds (cf. \cite[Exercise 13.27]{K}):
\begin{align}\label{eqn:modular}
\Xi^{(m)}_{r,s}(\frac{-1}{z})&=\sqrt{\frac{8}{(m+2)(m+3)}}\\
&\qquad \cdot\sum_{1\leq j\leq i\leq m+1} (-1)^{(r+s)(i+j)}
\sin\frac{\pi ri}{m+2} \sin\frac{\pi
sj}{m+3}\Xi^{(m)}_{i,j}(z).\nonumber
\end{align}

Let $\eta(z)=q^{1/24}\prod_{i=1}^{\infty}(1-q^i)$ be the Dedekind
$\eta$-function. The following transformation formula is
well-known (cf. \cite{Ap}):
\begin{equation*}
\eta(\frac{-1}{z}) = (-\sqrt{-1}z)^{1/2}\eta(z),
\end{equation*}
where we choose the branch of the square root function $x^{1/2}$
so that it is positive when $x>0$.

We review notations and some properties of the trace function in
\cite{DLM}. Let $g,h\in \Aut (M)$ be such that $gh=hg$. Let ${\cal
C}_1(g,h)$ be the space of $(g,h)$ $1$-point functions. Let $W$ be
a $g$-twisted $h$-stable $M$-module with conformal weight
$\lambda$. There is a linear isomorphism $\phi(h) : W\rightarrow
W$ such that
\begin{equation*}
\phi(h)Y_W(u,z)=Y_W(hu,z)\phi(h).
\end{equation*}

Define
\begin{equation*}
T_{W}(u,(g,h),z) = \tr_W u_{\wt(u)-1}\phi(h)q^{L(0)-1/20}
\end{equation*}
for homogeneous $u\in M$ and extend it for arbitrary $u \in M$
linearly. Note that the central charge of $M$ is $6/5$. Then
$T_{W}(\,\cdot\,,(g,h),z) \in {\cal C}_1(g,h)$ by \cite[Theorem
8.1]{DLM}. Let $F(\,\cdot\,,z)\in {\cal C}_1(g,h)$ and
$A=\left(\begin{array}{cc}a&b\\c&d\end{array}\right) \in
SL_2({\mathbb Z})$. Define $F|_A$ by
\begin{equation*}
F|_A(u,z)=(cz+d)^{-k}F\big(u,\frac{az+b}{cz+d}\big)
\end{equation*}
for $u\in M_{[k]}$ and extend it for arbitrary $u \in M$ linearly.
Then $F|_A\in {\cal C}_1(g^ah^c,g^bh^d)$ by \cite[Theorem
5.4]{DLM}. We denote $T_{W}({\bf 1},(g,h),z)$ by $T_W((g,h),z)$
for simplicity. Recall that the character $\ch W$ of $W$ is
defined to be $\tr_W q^{L(0)}$.

We want to determine the characters of the $20$ irreducible
$M^{\tau}$-modules $M(\epsilon)$, $W(\epsilon)$, $M_k^c$, $W_k^c$,
$M_T(\tau^i)(\epsilon)$, and $W_T(\tau^i)(\epsilon)$ for $\epsilon
= 0,1,2$ and $i=1,2$. We have shown in Theorem 2.1 that $\Aut (M)$
is generated by $\sigma$ and $\tau$. We shall consider the cases
where $g=1$ and $h=\tau$ or $g=\tau$ and $h=1$. We specify
$\phi(h)$ as follows. If $h=1$, we take $\phi(h)=1$. We shall deal
with the case $g=1$ and $h=\tau$ for $W=M$ or $W_k^0$. In such a
case we consider the same $\phi(\tau)$ as in Section 4. Thus if
$W=M$, we take $\phi(\tau)$ to be the automorphism $\tau$. If
$W=W_k^0$, we take $\phi(\tau)$ to be the linear isomorphism which
is naturally induced from the isometry $\tau$ of the lattice $(L,
\la\,\cdot\,,\,\cdot\,\ra)$.

Note that $T_W((g,1),z)=q^{-1/20}\ch W$. Note also that the
symmetry \eqref{sym2} induced by $\sigma$ implies
$T_{M(1)}((1,1),z)=T_{M(2)}((1,1),z)$. A similar assertion holds
for $W(1)$ and $W(2)$.

\begin{prop}
For $i=1,2$,
\begin{align*}
T_{M_T(\tau^i)}((\tau^i,1),z) &= \frac{\eta(z)}{\eta(z/3)}
(-\Xi^{(3)}_{2,1}-\Xi^{(3)}_{3,1}+\Xi^{(3)}_{3,3}),\\
T_{W_T(\tau^i)}((\tau^i,1),z) &= \frac{\eta(z)}{\eta(z/3)}
(\Xi^{(3)}_{1,1}+\Xi^{(3)}_{4,1}-\Xi^{(3)}_{4,3}).
\end{align*}
\end{prop}

\prf  Since $\ch V_L^{T_{\chi_j}}(\tau) = \ch S(\tau)$ for
$j=0,1,2$, we have
\begin{equation*}
q^{-1/12} \ch V_L^{T_{\chi_j}}(\tau)=\frac{\eta(z)}{\eta(z/3)}
\end{equation*}
by \eqref{chs}. Then \eqref{decomp1} and \eqref{decomp2} imply
that
\begin{align}
\frac{\eta(z)}{\eta(z/3)} &= T_{M_T(\tau)}((\tau,1),z)\cdot
(\Xi^{(3)}_{1,1}+\Xi^{(3)}_{4,1})+
T_{W_T(\tau)}((\tau,1),z)\cdot(\Xi^{(3)}_{2,1}+\Xi^{(3)}_{3,1}),
\label{chi0}\\
\frac{\eta(z)}{\eta(z/3)} &=
T_{M_T(\tau)}((\tau,1),z)\cdot\Xi^{(3)}_{4,3}+
T_{W_T(\tau)}((\tau,1),z)\cdot\Xi^{(3)}_{3,3}.\label{chi1}
\end{align}

Now consider
$\big(\Xi^{(3)}_{1,1}+\Xi^{(3)}_{4,1}\big)\Xi^{(3)}_{3,3}
-\big(\Xi^{(3)}_{2,1}+\Xi^{(3)}_{3,1}\big)\Xi^{(3)}_{4,3}$. Using
\eqref{eqn:modular} we can verify that it is invariant under the
action of $SL_2({\mathbb Z})$. Moreover, its $q$-expansion is
$1+0\cdot q+\cdots $. Thus
\begin{equation*}
(\Xi^{(3)}_{1,1}+\Xi^{(3)}_{4,1})\Xi^{(3)}_{3,3}
-(\Xi^{(3)}_{2,1}+\Xi^{(3)}_{3,1})\Xi^{(3)}_{4,3}=1.
\end{equation*}

Hence the assertions for $i=1$ follow from \eqref{chi0} and
\eqref{chi1}. The assertions for $i=2$ also hold by the symmetry
\eqref{sym2}. \quad \qed

\begin{thm}\label{theorem : character}
The characters of the $20$ irreducible $M^{\tau}$-modules
$M(\epsilon)$, $W(\epsilon)$, $M_k^c$, $W_k^c$,
$M_T(\tau^i)(\epsilon)$, and $W_T(\tau^i)(\epsilon)$ for $\epsilon
= 0,1,2$ and $i=1,2$ are given by the following formulas:

$(1)$ For $\epsilon=1,2$ we have
\begin{align*}
q^{-1/20}\ch M(0) & =
\frac{1}{3}\Big(\Xi^{(1)}_{1,1}\Xi^{(2)}_{1,1}+
\Xi^{(1)}_{2,1}\Xi^{(2)}_{3,1}
+2\frac{\eta(z)}{\eta(3z)}\Xi^{(3)}_{3,3}\Big),\\
q^{-1/20}\ch M(\epsilon) &= \frac{1}{3} \Big(
\Xi^{(1)}_{1,1}\Xi^{(2)}_{1,1}+
\Xi^{(1)}_{2,1}\Xi^{(2)}_{3,1}
-\frac{\eta(z)}{\eta(3z)}\Xi^{(3)}_{3,3}\Big),\\
q^{-1/20}\ch W(0) &= \frac{1}{3}\Big(
\Xi^{(1)}_{1,1}\Xi^{(2)}_{3,2}+
\Xi^{(1)}_{2,1}\Xi^{(2)}_{3,3}
-2\frac{\eta(z)}{\eta(3z)}\Xi^{(3)}_{4,3}\Big),\\
q^{-1/20}\ch W(\epsilon) &= \frac{1}{3} \Big(
\Xi^{(1)}_{1,1}\Xi^{(2)}_{3,2}+
\Xi^{(1)}_{2,1}\Xi^{(2)}_{3,3}+
\frac{\eta(z)}{\eta(3z)}\Xi^{(3)}_{4,3}\Big),\\
q^{-1/20}\ch M^c_k &= \Xi^{(1)}_{22}\Xi^{(2)}_{21},\\
q^{-1/20}\ch W^c_k &= \Xi^{(1)}_{22}\Xi^{(2)}_{22}.\\
\end{align*}

$(2)$ For $i=1,2$ we have
\begin{align*}
\left(\begin{array}{c}
q^{-1/20}\ch (M_T(\tau^i)(0))\\
q^{-1/20}\ch (M_T(\tau^i)(1))\\
q^{-1/20}\ch (M_T(\tau^i)(2))
\end{array}
\right) &\!=\! \frac{1}{3} \left(\begin{array}{ccc}
1 & 1 & 1\\
1 & \xi & \xi^2\\
1 & \xi^2 & \xi
\end{array}
\right) \left(\begin{array}{r}
T_{M_T(\tau^i)}((\tau^i,1),z)\\
e^{-11\pi\sqrt{-1}/90}T_{M_T(\tau^i)}((\tau^i,1),z+1)\\
e^{-22\pi\sqrt{-1}/90}T_{M_T(\tau^i)}((\tau^i,1),z+2)
\end{array}
\right),\\
\ &  \\
\left(\begin{array}{c}
q^{-1/20}\ch (W_T(\tau^i)(0))\\
q^{-1/20}\ch (W_T(\tau^i)(1))\\
q^{-1/20}\ch (W_T(\tau^i)(2))
\end{array}
\right) &\!=\! \frac{1}{3} \left(\begin{array}{ccc}
1 & 1 & 1\\
1 & \xi^2 & \xi\\
1 & \xi & \xi^2
\end{array}
\right) \left(\begin{array}{r}
T_{W_T(\tau^i)}((\tau^i,1),z)\\
e^{\pi\sqrt{-1}/90}T_{W_T(\tau^i)}((\tau^i,1),z+1)\\
e^{2\pi\sqrt{-1}/90}T_{W_T(\tau^i)}((\tau^i,1),z+2)
\end{array}
\right),
\end{align*}
where $\xi = \exp(2\pi\sqrt{-1}/3)$.
\end{thm}

\prf Since $M_T(\tau^i) = \oplus_{\epsilon=0}^2
M_T(\tau^i)(\epsilon)$ for $i=1,2$, we have
\begin{equation*}
T_{M_T(\tau^i)}((\tau^i,1),z) = \sum_{\epsilon=0}^2
T_{M_T(\tau^i)(\epsilon)}((1,1),z).
\end{equation*}

Replace $z$ with $z+k$, where $k=0,1,2$. Then
\begin{equation*}
T_{M_T(\tau^i)(\epsilon)}((1,1),z+k) = \tr_{M_T(\tau^i)(\epsilon)}
q^{L(0)-1/20} \exp(2\pi\sqrt{-1}k)^{L(0)-1/20}.
\end{equation*}

Note that $\exp(2\pi\sqrt{-1}k)^{L(0)-1/20} =
\exp(11\pi\sqrt{-1}k/90) \xi^{2k\epsilon}$ on
$M_T(\tau^i)(\epsilon)$, since the eigenvalues for $L(0)$ on
$M_T(\tau^i)(\epsilon)$ are of the form $1/9+2\epsilon/3+n$ with
$n \in \Z_{\ge 0}$. Thus
\begin{equation*}
T_{M_T(\tau^i)}((\tau^i,1),z+k) = \exp(11\pi\sqrt{-1}k/90)
\sum_{\epsilon=0}^2 \xi^{2k\epsilon}
T_{M_T(\tau^i)(\epsilon)}((1,1),z).
\end{equation*}

We can solve these equations for $k=0,1,2$ with respect to
$T_{M_T(\tau^i)(\epsilon)}((1,1),z)$, $\epsilon=0,1,2$, and obtain
the expressions of $T_{M_T(\tau^i)(\epsilon)}((1,1),z) = q^{-1/20}
\ch (M_T(\tau^i)(\epsilon))$ in the theorem.

Similarly, $W_T(\tau^i) = \oplus_{\epsilon=0}^2
W_T(\tau^i)(\epsilon)$ and the eigenvalues for $L(0)$ on
$W_T(\tau^i)(\epsilon)$ are of the form $2/45+\epsilon/3+n$ with
$n \in \Z_{\ge 0}$. Hence $\exp(2\pi\sqrt{-1}k)^{L(0)-1/20}\! =\!
\exp(-\pi\sqrt{-1}k/90) \xi^{k\epsilon}$ on
$W_T(\tau^i)(\epsilon)$ and we obtain the expressions of
$q^{-1/20} \ch (W_T(\tau^i)(\epsilon)), \epsilon=0,1,2$.

It is proved in \cite{LY} that $M=M_k^0$ is a rational vertex
operator algebra. Moreover, there are exactly two inequivalent
irreducible $\tau$-stable $M$-modules, namely, $M$ and $W_k^0$ by
Lemma 4.1. Since $M_T(\tau)$ and $W_T(\tau)$ are two inequivalent
irreducible $\tau$-twisted $M$-modules, we have $\dim {\cal
C}_1(\tau,1)=\dim {\cal C}_1(1,\tau)=2$ and
\begin{equation*}
\{T_{M_T(\tau)}(\,\cdot\,,(\tau,1),z),
T_{W_T(\tau)}(\,\cdot\,,(\tau,1),z)\}
\end{equation*}
is a basis of ${\cal
C}_1(\tau,1)$ by \cite[Theorems 5.4 and 10.1]{DLM}. Now
$T_M(\,\cdot\,,(1,\tau),z) |_{S}\in {\cal C}_{1}(\tau,1)$ for
$S=\left(\begin{array}{rr}0&-1\\1&0\end{array}\right)$ by
\cite[Theorems 5.4 and 8.1]{DLM}. Thus,
\begin{equation*}
T_M((1,\tau),z)  = \alpha T_{M_T(\tau)}((\tau,1),\frac{-1}{z})+
\beta T_{W_T(\tau)}((\tau,1),\frac{-1}{z})
\end{equation*}
for some $\alpha,\beta\in {\mathbb C}$.

>From \eqref{eqn:modular} and Proposition 6.1 it follows that
\begin{align*}
\frac{\eta(3z)}{\eta(z)}T_{M_T(\tau)}((\tau,1),\frac{-1}{z}) &=
\frac{2\sin(\frac{\pi}{5})}{\sqrt{5}}\Xi^{(3)}_{3,3}
-\frac{2\sin(\frac{2\pi}{5})}{\sqrt{5}}\Xi^{(3)}_{4,3},\\
\frac{\eta(3z)}{\eta(z)}T_{W_T(\tau)}((\tau,1),\frac{-1}{z})
&= \frac{2\sin(\frac{2\pi}{5})}{\sqrt{5}}\Xi^{(3)}_{3,3}
+\frac{2\sin(\frac{\pi}{5})}{\sqrt{5}}\Xi^{(3)}_{4,3}.
\end{align*}

Thus
\begin{equation*}
T_M((1,\tau),z) = q^{-1/20}\Big((\alpha
\frac{2\sin(\frac{\pi}{5})}{\sqrt{5}}+
\beta\frac{2\sin(\frac{2\pi}{5})}{\sqrt{5}}) +(-\alpha
\frac{2\sin(\frac{2\pi}{5})}{\sqrt{5}}+
\beta\frac{2\sin(\frac{\pi}{5})}{\sqrt{5}})q^{3/5}+\cdots\Big).
\end{equation*}

Furthermore, we see that $T_M((1,\tau),z)=q^{-1/20}( 1+0\cdot
q^{3/5}+\cdots)$ by a direct computation. Hence
$\alpha=2\sin(\frac{\pi}{5})/\sqrt{5}$ and
$\beta=2\sin(\frac{2\pi}{5})/\sqrt{5}$. Therefore,
\begin{equation*}
T_M((1,\tau),z) = \frac{\eta(z)}{\eta(3z)}\Xi^{(3)}_{3,3}.
\end{equation*}

Note that
\begin{equation*} T_M((1,1),z)=q^{-1/20}\ch M =
\Xi^{(1)}_{1,1}\Xi^{(2)}_{1,1}+\Xi^{(1)}_{2,1}\Xi^{(2)}_{3,1}.
\end{equation*}

Now $M=M(0)\oplus M(1)\oplus M(2)$ and $T_{M(1)}((1,1),z) =
T_{M(2)}((1,1),z)$ by the symmetry \eqref{sym2}. Then
\begin{align*}
T_M((1,1),z) &= T_{M(0)}((1,1),z)+ T_{M(1)}((1,1),z)+
T_{M(2)}((1,1),z)\\
&= T_{M(0)}((1,1),z)+2T_{M(1)}((1,1),z)
\end{align*}
and
\begin{align*}
T_M((1,\tau),z) &= T_{M(0)}((1,1),z)+ \xi T_{M(1)}((1,1),z)+
\xi^2 T_{M(2)}((1,1),z)\\
&= T_{M(0)}((1,1),z)-T_{M(1)}((1,1),z)
\end{align*}
by the definition of trace functions. Thus $q^{-1/20} \ch
M(\epsilon)=T_{M(\epsilon)}((1,1),z)$ can be expressed as
\begin{align*}
q^{-1/20}\ch M(0)
&= \frac{1}{3}\Big(T_M((1,1),z)+2T_M((1,\tau),z)\Big)\\
&= \frac{1}{3}\Big( \Xi^{(1)}_{1,1}\Xi^{(2)}_{1,1}+
\Xi^{(1)}_{2,1}\Xi^{(2)}_{3,1}+
2\frac{\eta(z)}{\eta(3z)}\Xi^{(3)}_{3,3}\Big),\\
q^{-1/20}\ch M(\epsilon)
&= \frac{1}{3}\Big( T_M((1,1),z)-T_M((1,\tau),z)\Big)\\
&= \frac{1}{3}\Big( \Xi^{(1)}_{1,1}\Xi^{(2)}_{1,1}+
\Xi^{(1)}_{2,1}\Xi^{(2)}_{3,1}-
\frac{\eta(z)}{\eta(3z)}\Xi^{(3)}_{3,3}\Big)
\end{align*}
for $\epsilon = 1,2$. The computations for $W(\epsilon),
\epsilon=0,1,2$ are similar.

Since $M_k^i$, $i=a,b,c$ are equivalent irreducible
$M^{\tau}$-modules by Lemma 4.1, we have $q^{-1/20}\ch M_k^c =
q^{-1/20}\ch M_k^a = \Xi^{(1)}_{2,2} \Xi^{(2)}_{2,1}$. Likewise,
$q^{-1/20}\ch W_k^c = q^{-1/20}\ch W_k^a = \Xi^{(1)}_{2,2}
\Xi^{(2)}_{2,2}$. The proof is complete. \quad \qed

\medskip
We now disccuss the relation between the characters computed here
and those of modules for a $W$-algebra computed in \cite{FKW}. We
use the notation of \cite{FKW} without any comments. We refer to
their results in the case that $\bar{\g}$ is the simple finite
dimensional Lie algebra over ${\mathbb C}$ of type $A_2$ and
$(p,p^{\prime})=(6,5)$. In this case, we have
\begin{align*}
P^{p-h^{\vee}}_{+} =P^{3}_{+} &= \{\sum_{i=0}^{2}a_i
\Lambda_i\,|\, 0\leq a_i\in {\mathbb Z} \mbox{ and }
\sum_{i=0}^{2}a_i=3\},\\
{P^{\vee}_{+}}^{p^{\prime}-h}={P^{\vee}_{+}}^{2} &=
\{\sum_{i=0}^{2}b_i\Lambda_i^{\vee}\, |\, 0\leq b_i\in {\mathbb Z}
\mbox{ and } \sum_{i=0}^{2}b_i=2\}.
\end{align*}

It can be easily shown that $\widetilde{W}_{+}=\langle g\rangle$
is the cyclic group of order $3$ such that
$g(\Lambda_0)=\Lambda_1, g(\Lambda_1)=\Lambda_2$, and
$g(\Lambda_2)=\Lambda_0$. The cardinality of $I_{p,p^{\prime}}=
(P_{+}^{3}\times {P^{\vee}_{+}}^{2})/\widetilde{W}_{+}$ is equal
to $20$.

For $\lambda\in P^{3}_{+},
\mu\in {P^{\vee}_{+}}^{2}$, define
\begin{equation*}
\varphi_{\lambda,\mu}(z) =  \eta(z)^{-2}\sum_{w\in W}\epsilon(w)
q^{\frac{1}{2pp^{\prime}}
|p^{\prime}w(\lambda+\rho)-p(\mu+\rho^{\vee})|^2}.
\end{equation*}

The vector space spanned by $\varphi_{\lambda,\mu}(z),
(\lambda,\mu)\in I_{p,p^{\prime}}$ is invariant under the action
of $SL_2({\mathbb Z})$ and the transformation formula
\begin{equation*}
\varphi_{\lambda, \lambda^{\prime}}(\frac{-1}{z})=
\sum_{(\mu,\mu^{\prime})\in I_{p,p^{\prime}}} S_{(\lambda,
\lambda^{\prime}),(\mu,\mu^{\prime})}
\varphi_{\mu,\mu^{\prime}}(z)
\end{equation*}
is given by \cite[(4.2.2)]{FKW}. Define ${\cal
F}_1=\{\varphi_{\lambda,\mu}(z) \ |\ (\lambda,\mu) \in
I_{p,p^{\prime}}\}$. In \cite[Section 3]{FKW}, it is shown that
each $\varphi_{\lambda,\mu}(z)\in {\cal F}_1$ is the character of
a module for the $W$-algebra associated to $\bar{\g}$ and
$(p,p^{\prime})$ which is conjectured to be irreducible.

We denote by ${\cal F}_2$ the set of characters of all irreducible
$M^{\tau}$-modules computed in Theorem \ref{theorem : character}.
For any $m$, there is a congruence subgroup $\Gamma_m$ such that
each $\Xi^{(m)}_{r,s}$ is a modular form for $\Gamma_m$ (cf.
\cite[(6.11)]{W}). Then there is a congruence subgroup $\Gamma$
such that each character in ${\cal F}_2$ is invariant under the
action of $\Gamma$.
The following transformation formulas hold by the formula
\eqref{eqn:modular}:
\begin{equation*}
\left(\begin{array}{c}
T_{M}((1,1),\frac{-1}{z})\\
T_{W^0_k}((1,1),\frac{-1}{z})\\
T_{M^c_k}((1,1),\frac{-1}{z})\\
T_{W^c_k}((1,1),\frac{-1}{z})
\end{array}
\right)= \left(
\begin{array}{rrrrrr}
\frac{\sin(\frac{\pi}{5})}{\sqrt{5}} &
\frac{\sin(\frac{2\pi}{5})}{\sqrt{5}} &
\frac{3\sin(\frac{\pi}{5})}{\sqrt{5}} &
\frac{3\sin(\frac{2\pi}{5})}{\sqrt{5}}\\
\frac{\sin(\frac{2\pi}{5})}{\sqrt{5}} &
-\frac{\sin(\frac{\pi}{5})}{\sqrt{5}} &
\frac{3\sin(\frac{2\pi}{5})}{\sqrt{5}} &
-\frac{3\sin(\frac{\pi}{5})}{\sqrt{5}}\\
\frac{\sin(\frac{\pi}{5})}{\sqrt{5}} &
\frac{\sin(\frac{2\pi}{5})}{\sqrt{5}} &
-\frac{\sin(\frac{\pi}{5})}{\sqrt{5}} &
-\frac{\sin(\frac{2\pi}{5})}{\sqrt{5}} \\
\frac{\sin(\frac{2\pi}{5})}{\sqrt{5}} &
-\frac{\sin(\frac{\pi}{5})}{\sqrt{5}} &
-\frac{\sin(\frac{2\pi}{5})}{\sqrt{5}} &
\frac{\sin(\frac{\pi}{5})}{\sqrt{5}}
\end{array}
\right) \left(\begin{array}{c}
T_{M}((1,1),z)\\
T_{W^0_k}((1,1),z)\\
T_{M^c_k}((1,1),z)\\
T_{W^c_k}((1,1),z)
\end{array}
\right)
\end{equation*}
and
\begin{equation*}
\left(\begin{array}{c}
T_{M}((1,\tau),\frac{-1}{z})\\
T_{W^0_k}((1,\tau),\frac{-1}{z})
\end{array}\right)
=
\left(\begin{array}{rr}
\frac{2\sin (\frac{\pi}{5})}{\sqrt{5}}&
\frac{2\sin (\frac{2\pi}{5})}{\sqrt{5}}\\
\frac{2\sin (\frac{2\pi}{5})}{\sqrt{5}} &
\frac{-2\sin (\frac{\pi}{5})}{\sqrt{5}}
\end{array}\right)
\left(\begin{array}{c}
T_{M_T(\tau)}((\tau,1),z)\\
T_{W_T(\tau)}((\tau,1),z)
\end{array}\right).
\end{equation*}

Thus we have the transformation formulas for elements of ${\cal
F}_2$. Comparing the $q$-expansions and the coefficients of
transformation formulas of elements in ${\cal F}_1$ and ${\cal
F}_2$, it can be shown that ${\cal F}_1={\cal F}_2$ using Lemma
1.7.1 in \cite{KW}. In particular,
$\varphi_{3\Lambda_0,2\Lambda_0^{\vee}}(z) =q^{-1/20}\ch M^{\tau}$
holds.

\newpage

\appendix

\section{$\bv^{12}$, $J(-1)\bv^{12}$, $J(-2)\bv^{12}$, and
$J(-1)^2\bv^{12}$}

{\footnotesize
\begin{align*}
\bv^{12} &=
-(5877264800/3501)L(-12)\1+(3404072000/3501)L(-10)L(-2)\1\\
& \quad-(2653990000/3501)L(-9)L(-3)\1-(266376800/3501)L(-8)L(-4)\1\\
& \quad +(282988000/1167)L(-8)L(-2)^2\1-(23744800/1167)L(-7)L(-5)\1\\
& \quad -(30824000/1167)L(-7)L(-3)L(-2)\1+(1242377600/1167)L(-6)^2\1\\
& \quad -(61947200/3501)L(-6)L(-4)L(-2)\1-(1313806000/1167)L(-6)L(-3)^2\1 \\
& \quad -(45496000/1167)L(-6)L(-2)^3\1-(3046768400/3501)L(-5)^2L(-2)\1\\
& \quad +(299424800/1167)L(-5)L(-4)L(-3)\1+(2347094000/3501)L(-5)L(-3)L(-2)^2\1\\
& \quad -(17280400/1167)L(-4)^3\1-(2036373200/3501)L(-4)^2L(-2)^2\1\\
& \quad +(82996000/3501)L(-4)L(-3)^2L(-2)\1+(1074512000/3501)L(-4)L(-2)^4\1\\
& \quad +(511628125/3501)L(-3)^4\1-(418850000/3501)L(-3)^2L(-2)^3\1\\
& \quad -(59680000/3501)L(-2)^6\1-(505200/389)L(-6)J(-3)^2\1\\
& \quad +(3380480/1167)L(-4)L(-2)J(-3)^2\1+1150L(-3)^2J(-3)^2\1 \\
& \quad -(184400/1167)L(-2)^3J(-3)^2\1+(3788680/1167)L(-5)J(-4)J(-3)\1 \\
& \quad -(8788400/3501)L(-3)L(-2)J(-4)J(-3)\1-(12761440/3501)L(-4)J(-5)J(-3)\1 \\
& \quad -(5727500/10503)L(-4)J(-4)^2\1+(352400/389)L(-2)^2J(-5)J(-3)\1 \\
& \quad +(5727500/10503)L(-2)^2J(-4)^2\1+(1593900/389)L(-3)J(-6)J(-3)\1 \\
& \quad +(12935800/10503)L(-3)J(-5)J(-4)\1+(4108000/3501)L(-2)J(-7)J(-3)\1 \\
& \quad -(2811800/1167)L(-2)J(-6)J(-4)\1-(3131600/10503)L(-2)J(-5)^2\1 \\
& \quad -(14904160/3501)J(-9)J(-3)\1+(32677600/10503)J(-8)J(-4)\1 \\
& \quad +(9423200/10503)J(-7)J(-5)\1+(2432375/1167)J(-6)^2\1\\
& \quad +J(-3)^4\1.
\end{align*}
}

{\footnotesize
\begin{align*}
J(-1)\bv^{12} & =
(47528/389)L(-4)J(-3)^3\1-(53552200/1167)J(-7)J(-3)^2\1\\
& \quad -(14322122880/389)L(-10)J(-3)\1-(7313862400/389)L(-8)L(-2)J(-3)\1\\
& \quad -(2263268800/389)L(-7)L(-3)J(-3)\1 +(7140323840/1167)L(-6)L(-4)J(-3)\1\\
& \quad -(4870066240/389)L(-5)^2J(-3)\1 -(41271174880/3501)L(-9)J(-4)\1\\
& \quad +(87811701120/389)J(-13)\1 -(65647722400/389)L(-2)J(-11)\1\\
& \quad +(195884000/1167)L(-8)J(-5)\1 +(18292448200/389)L(-3)J(-10)\1\\
& \quad -(2704504400/389)L(-7)J(-6)\1 -(8342231040/389)L(-4)J(-9)\1\\
& \quad -(2134787200/389)L(-6)J(-7)\1+(17270275360/1167)L(-5)J(-8)\1\\
& \quad -(926640/389)L(-2)^2J(-3)^3\1 +(24986000/1167)L(-2)J(-5)J(-3)^2\1\\
& \quad -(2833833600/389)L(-6)L(-2)^2J(-3)\1 +(1692496000/389)L(-5)L(-3)L(-2)J(-3)\1\\
& \quad -(6705813920/1167)L(-4)^2L(-2)J(-3)\1 -(10899200/9)L(-7)L(-2)J(-4)\1\\
& \quad +(20147275200/389)L(-2)^2J(-9)\1 +(32842950400/3501)L(-6)L(-2)J(-5)\1\\
& \quad -(8472651200/389)L(-3)L(-2)J(-8)\1 -(10511649200/1167)L(-5)L(-2)J(-6)\1\\
& \quad +(12944796800/1167)L(-4)L(-2)J(-7)\1 -(9963200/389)J(-5)^2J(-3)\1\\
& \quad +(1607444800/389)L(-5)L(-3)J(-5)\1 -(1408915040/1167)L(-4)^2J(-5)\1\\
& \quad +(1271140/1167)L(-3)J(-4)J(-3)^2\1 -(2312728600/1167)L(-4)L(-3)^2J(-3)\1\\
& \quad -(38446491200/3501)L(-6)L(-3)J(-4)\1 -(1002365000/1167)L(-3)^2J(-7)\1\\
& \quad +(462527000/389)L(-4)L(-3)J(-6)\1 +(14283100/389)J(-6)J(-4)J(-3)\1\\
& \quad +(15075473920/3501)L(-5)L(-4)J(-4)\1 -(11661041600/3501)L(-4)L(-2)^2J(-5)\1\\
& \quad -(7223710000/3501)L(-3)^2L(-2)J(-5)\1 +(2639390000/389)L(-3)L(-2)^2J(-6)\1\\
& \quad -(30846400/389)L(-5)L(-2)^2J(-4)\1 +(5365349600/3501)L(-4)L(-3)L(-2)J(-4)\1\\
& \quad -(34590712000/3501)L(-2)^3J(-7)\1 +(6230282500/3501)L(-3)^3J(-4)\1\\
& \quad +(1547000/389)J(-5)J(-4)^2\1 -(4562948000/3501)L(-3)L(-2)^3J(-4)\1\\
& \quad -(10829000/1167)L(-2)J(-4)^2J(-3)\1 +(4340336000/3501)L(-2)^4J(-5)\1\\
& \quad +(1919403200/1167)L(-4)L(-2)^3J(-3)\1-(99658000/389)L(-3)^2L(-2)^2J(-3)\1\\
& \quad -(89856000/389)L(-2)^5J(-3)\1.
\end{align*}
}

{\footnotesize
\begin{align*}
J(-2)\bv^{12} &=
-4272L(-5)J(-3)^3\1-(21069744/389)J(-8)J(-3)^2\1\\
& \quad -(14150438080/1167)L(-11)J(-3)\1 -(3639849600/389)L(-9)L(-2)J(-3)\1\\
& \quad -(9699222400/1167)L(-8)L(-3)J(-3)\1 +(2157139840/1167)L(-7)L(-4)J(-3)\1\\
& \quad -(5925448960/1167)L(-6)L(-5)J(-3)\1 +(3006435200/389)L(-10)J(-4)\1\\
& \quad +(325064548960/389)J(-14)\1 -(176823168000/389)L(-2)J(-12)\1\\
& \quad +(10174691200/3501)L(-9)J(-5)\1 +(38988751200/389)L(-3)J(-11)\1\\
& \quad -(4612321600/1167)L(-8)J(-6)\1 -(33023056960/389)L(-4)J(-10)\1\\
& \quad -(13371577600/1167)L(-7)J(-7)\1 +(54711326720/1167)L(-5)J(-9)\1\\
& \quad +(4368409600/1167)L(-6)J(-8)\1 -960L(-3)L(-2)J(-3)^3\1\\
& \quad +(5080480/389)L(-2)J(-6)J(-3)^2\1 -(1269523200/389)L(-7)L(-2)^2J(-3)\1\\
& \quad +(1626342400/1167)L(-6)L(-3)L(-2)J(-3)\1 -(3954100480/1167)L(-5)L(-4)L(-2)J(-3)\1 \\
& \quad +(6597673600/1167)L(-8)L(-2)J(-4)\1+(382495595200/3501)L(-2)^2J(-10)\1 \\
& \quad +(23379344000/3501)L(-7)L(-2)J(-5)\1-(41865472000/1167)L(-3)L(-2)J(-9)\1 \\
& \quad +(5662851200/1167)L(-6)L(-2)J(-6)\1+(41335582720/1167)L(-4)L(-2)J(-8)\1 \\
& \quad -(66974297600/3501)L(-5)L(-2)J(-7)\1+7760L(-3)J(-5)J(-3)^2\1 \\
& \quad -(13118000/1167)J(-6)J(-5)J(-3)\1-(3036691200/389)L(-6)L(-3)J(-5)\1 \\
& \quad +(2541514240/1167)L(-5)L(-4)J(-5)\1-(4489884400/1167)L(-5)L(-3)^2J(-3)\1\\
& \quad +(48898000/389)L(-5)L(-3)J(-6)\1+(524720000/389)L(-4)^2L(-3)J(-3)\1 \\
& \quad -(478727200/389)L(-4)^2J(-6)\1-(315678400/3501)L(-7)L(-3)J(-4)\1 \\
& \quad -(5656762000/1167)L(-3)^2J(-8)\1+(4388915200/1167)L(-4)L(-3)J(-7)\1 \\
& \quad +(5080480/1167)L(-4)J(-4)J(-3)^2\1+(7809478400/1167)L(-6)L(-4)J(-4)\1 \\
& \quad +(117493120/3501)J(-7)J(-4)J(-3)\1-(6924715520/3501)L(-5)^2J(-4)\1 \\
& \quad +(7972739200/3501)L(-5)L(-2)^2J(-5)\1+(726208000/3501)L(-4)L(-3)L(-2)J(-5)\1 \\
& \quad +(15160000/3501)J(-5)^2J(-4)\1-(9229774400/1167)L(-4)L(-2)^2J(-6)\1 \\
& \quad -(1273060000/1167)L(-3)^2L(-2)J(-6)\1-(5021408000/3501)L(-6)L(-2)^2J(-4)\1 \\
& \quad +(4736835200/1167)L(-5)L(-3)L(-2)J(-4)\1-(4697646400/3501)L(-4)^2L(-2)J(-4)\1 \\
& \quad -(28325800/3501)J(-6)J(-4)^2\1+(10330016000/1167)L(-3)L(-2)^2J(-7)\1 \\
& \quad +(6184910000/3501)L(-3)^3J(-5)\1-(2988476000/3501)L(-4)L(-3)^2J(-4)\1 \\
& \quad +(2298688000/1167)L(-4)L(-2)^3J(-4)\1-(59886716800/3501)L(-2)^3J(-8)\1 \\
& \quad -(1320284000/3501)L(-3)^2L(-2)^2J(-4)\1+(22910000/10503)L(-2)J(-4)^3\1 \\
& \quad -(3979216000/3501)L(-3)L(-2)^3J(-5)\1-(122435200/389)L(-5)L(-2)^3J(-3)\1 \\
& \quad -(977670400/1167)L(-4)L(-3)L(-2)^2J(-3)\1-(22467200/3501)L(-2)J(-5)J(-4)J(-3)\1 \\
& \quad +(58888000/1167)L(-3)^3L(-2)J(-3)\1-(17576800/3501)L(-3)J(-4)^2J(-3)\1 \\
& \quad +(29504000/389)L(-3)L(-2)^4J(-3)\1+(2281792000/1167)L(-2)^4J(-6)\1 \\
& \quad -(368800/389)L(-2)^2J(-4)J(-3)^2\1-(238720000/1167)L(-2)^5J(-4)\1.
\end{align*}
}

{\footnotesize
\begin{align*}
J(-1)^2\bv^{12} &=
(28587850894720/389)L(-14)\1+(40679435680000/1167)L(-12)L(-2)\1 \\
& \quad -(20370766707200/389)L(-11)L(-3)\1-(29040708661120/389)L(-10)L(-4)\1 \\
& \quad -(1357372140800/389)L(-10)L(-2)^2\1-(120978369778240/1167)L(-9)L(-5)\1 \\
& \quad +(15046999864000/1167)L(-9)L(-3)L(-2)\1-(120139236131200/1167)L(-8)L(-6)\1 \\
& \quad +(7353135836800/1167)L(-8)L(-4)L(-2)\1+(5914869272000/389)L(-8)L(-3)^2\1 \\
& \quad -(9027652192000/1167)L(-8)L(-2)^3\1-(19757556187200/389)L(-7)^2\1 \\
& \quad +(10357377908800/389)L(-7)L(-5)L(-2)\1+(6212435174400/389)L(-7)L(-4)L(-3)\1 \\
& \quad -(3066391744000/389)L(-7)L(-3)L(-2)^2\1-(34866323814400/1167)L(-6)^2L(-2)\1 \\
& \quad -(1360052761600/389)L(-6)L(-5)L(-3)\1-(3455809144320/389)L(-6)L(-4)^2\1 \\
& \quad +(8114060115200/1167)L(-6)L(-4)L(-2)^2\1+(2356317080000/1167)L(-6)L(-3)^2L(-2)\1 \\
& \quad -(4200302912000/1167)L(-6)L(-2)^4\1+(2046779720960/389)L(-5)^2L(-4)\1 \\
& \quad +(5012264899200/389)L(-5)^2L(-2)^2\1+(5606971697600/1167)L(-5)L(-4)L(-3)L(-2)\1 \\
& \quad +(4546296703000/1167)L(-5)L(-3)^3\1-(3986231288000/1167)L(-5)L(-3)L(-2)^3\1 \\
& \quad -(824891421120/389)L(-4)^3L(-2)\1+(129922182000/389)L(-4)^2L(-3)^2\1 \\
& \quad +(9190279446400/1167)L(-4)^2L(-2)^3\1-(3417631724000/1167)L(-4)L(-3)^2L(-2)^2\1 \\
& \quad -(1854416512000/1167)L(-4)L(-2)^5\1-(339474200000/1167)L(-3)^4L(-2)\1 \\
& \quad +(472407520000/1167)L(-3)^2L(-2)^4\1+(21565440000/389)L(-2)^7\1 \\
& \quad -(33906046720/389)L(-8)J(-3)^2\1-(38547928640/389)L(-6)L(-2)J(-3)^2\1 \\
& \quad +(8889576280/389)L(-5)L(-3)J(-3)^2\1-(52680368/389)L(-4)^2J(-3)^2\1 \\
& \quad +(1681515680/389)L(-4)L(-2)^2J(-3)^2\1-(4900781600/389)L(-3)^2L(-2)J(-3)^2\1 \\
& \quad -(680659200/389)L(-2)^4J(-3)^2\1-(21316634560/1167)L(-7)J(-4)J(-3)\1 \\
& \quad +(15456968800/389)L(-5)L(-2)J(-4)J(-3)\1-(57407779520/1167)L(-4)L(-3)J(-4)J(-3)\1 \\
& \quad +(769371200/389)L(-3)L(-2)^2J(-4)J(-3)\1+(82018834560/389)L(-6)J(-5)J(-3)\1 \\
& \quad -(318755320000/3501)L(-6)J(-4)^2\1-(62232722240/1167)L(-4)L(-2)J(-5)J(-3)\1 \\
& \quad +(59657182000/3501)L(-4)L(-2)J(-4)^2\1+(4384283800/1167)L(-3)^2J(-5)J(-3)\1 \\
& \quad +(28313585300/1167)L(-3)^2J(-4)^2\1+(14719931200/1167)L(-2)^3J(-5)J(-3)\1 \\
& \quad -(15017860000/3501)L(-2)^3J(-4)^2\1-(102815580920/389)L(-5)J(-6)J(-3)\1 \\
& \quad +(214806972640/3501)L(-5)J(-5)J(-4)\1+(20784972000/389)L(-3)L(-2)J(-6)J(-3)\1 \\
& \quad -(133605586400/3501)L(-3)L(-2)J(-5)J(-4)\1+(243575438080/1167)L(-4)J(-7)J(-3)\1 \\
& \quad +(7292932400/389)L(-4)J(-6)J(-4)\1-(12891781760/389)L(-4)J(-5)^2\1 \\
& \quad -(49983377600/389)L(-2)^2J(-7)J(-3)\1+(10825750000/389)L(-2)^2J(-6)J(-4)\1 \\
& \quad -(13957486400/3501)L(-2)^2J(-5)^2\1-(173848522640/1167)L(-3)J(-8)J(-3)\1 \\
& \quad -(65060216000/1167)L(-3)J(-7)J(-4)\1+(25622862200/389)L(-3)J(-6)J(-5)\1 \\
& \quad +(174271514560/389)L(-2)J(-9)J(-3)\1-(232573421600/3501)L(-2)J(-8)J(-4)\1 \\
& \quad +(392430209600/3501)L(-2)J(-7)J(-5)\1-(31534947600/389)L(-2)J(-6)J(-6)\1 \\
& \quad -(5559840/389)L(-2)J(-3)^4\1-(291151720080/389)J(-11)J(-3)\1 \\
& \quad +(257458099600/1167)J(-10)J(-4)\1-(140099797760/389)J(-9)J(-5)\1 \\
& \quad +(83988236280/389)J(-8)J(-6)\1-(44378890400/389)J(-7)J(-7)\1 \\
& \quad +(22538776/389)J(-5)J(-3)^3\1-(26131300/1167)J(-4)^2J(-3)^2\1.
\end{align*}
}

\section{The images of four vectors in $A(\W)$}

For simplicity of notation we omit the symbol $\ast$ for
multiplication in $A(\W)$.

\begin{align*}
[\bv^{12}] &=
-(59680000/3501)[\omega]^6+(156040000/3501)[\omega]^5
-(115878400/3501)[\omega]^4 \\
& \quad +\big(-(184400/1167)[J]^2+32328400/3501\big)[\omega]^3 \\
& \quad +\big((536500/1167)[J]^2-3155968/3501\big)[\omega]^2 \\
& \quad +\big(-(87812/389)[J]^2+93184/3501\big)[\omega] \\
& \quad +[J]^4+(75776/3501)[J]^2.\\
& \  \\
& \  \\
[J(-1)\bv^{12}] &=
-(89856000/389)[J][\omega]^5+(228945600/389)[J][\omega]^4 \\
& \quad -(555607520/1167)[J][\omega]^3 \\
& \quad +\big(-(926640/389)[J]^3+(57790304/389)[J]\big)[\omega]^2\\
& \quad +\big((1637064/389)[J]^3-(19542016/1167)[J]\big)[\omega] \\
& \quad -(668408/389)[J]^3+(186368/389)[J].\\
& \  \\
& \  \\
[J(-2)\bv^{12}] &=
(179712000/389)[J][\omega]^5-(457891200/389)[J][\omega]^4 \\
& \quad +(1111215040/1167)[J][\omega]^3 \\
& \quad +\big((1853280/389)[J]^3-(115580608/389)[J]\big)[\omega]^2 \\
& \quad +\big(-(3274128/389)[J]^3+(39084032/1167)[J]\big)[\omega] \\
& \quad +(1336816/389)[J]^3-(372736/389)[J].\\
& \  \\
& \  \\
[J(-1)^2\bv^{12}] &=
(21565440000/389)[\omega]^7+(513849856000/1167)[\omega]^6 \\
& \quad -(552497504000/389)[\omega]^5 \\
& \quad +\big(-(680659200/389)[J]^2+1285515063040/1167\big)[\omega]^4 \\
& \quad +\big((3994427840/389)[J]^2-121501591744/389\big)[\omega]^3 \\
& \quad +\big(-(8220864912/389)[J]^2+36103315456/1167\big)[\omega]^2 \\
& \quad +\big(-(5559840/389)[J]^4+(3836073072/389)[J]^2
-363417600/389\big)[\omega] \\
& \quad -(9879324/389)[J]^4-(355536896/389)[J]^2.
\end{align*}


\newpage
\begin{thebibliography}{99}
\bibitem{abd}
T. Abe, G. Buhl and C. Dong,
 Rationality, regularity and
$C_2$-cofiniteness, math.QA/0204021.

\bibitem{Ap}
T. M. Apostol,
Modular functions and Dirichlet series in number theory,
Second edition, Graduate Texts in Mathematics {\bfseries 41},
Springer-Verlag, New York, 1990.

\bibitem{BMP}
P. Bouwknegt, J. McCarthy and K. Pilch,
{\itshape The $\mathcal{W}_3$ Algebra}, Lecture Notes in Physics,
{\bfseries m42}, Springer, Berlin 1996.

\bibitem{buhl} G. Buhl, A spanning set for VOA modules,
{\itshape  J. Algebra} {\bfseries 254} (2002), 125--151.

\bibitem{D1}
C. Dong, Vertex algebras associated with even lattices, {\itshape J.
Algebra} {\bfseries 161} (1993), 245--265.

\bibitem{D2}
C. Dong, Twisted modules for vertex algebras associated with even
lattices, {\itshape J. Algebra} {\bfseries 165} (1994), 91--112.

\bibitem{DL}
C. Dong and J. Lepowsky, The algebraic structure of
relative twisted vertex operators, {\itshape J. Pure and Applied
Algebra} {\bfseries 110}(1996), 259--295.

\bibitem{dlmnew}
C. Dong, H. Li and  G. Mason, Twisted representations of vertex operator
algebras, {\itshape Math. Ann.} {\bfseries 310} (1998), 571--600.

\bibitem{DLM}
C. Dong, H. Li and  G. Mason, Modular-invariance of trace
functions in orbifold theory and generalized moonshine, {\itshape
Comm. Math. Phys.} {\bfseries 214} (2000), 1--56.

\bibitem{DLMN}
C. Dong, H. Li, G. Mason and S. P. Norton,
Associative subalgebras of the Griess algebra and related topics,
in: {\itshape Proc. of the Conference on the Monster and Lie algebras at
The Ohio State University, May 1996,}  ed. by J. Ferrar and K. Harada,
Walter de Gruyter, Berlin-New York, 1998, 27--42.

\bibitem{DM}
C. Dong and G. Mason, On quantum Galois theory,
{\itshape Duke Math. J.} {\bfseries 86} (1997), 305--321.

\bibitem{DN1}
C. Dong and K. Nagatomo, Classification of irreducible modules for
the vertex operator algebra $M(1)^+$, {\itshape J. Algebra}
{\bfseries 216} (1999), 384--404.

\bibitem{DY}
C. Dong and G. Yamskulna, Vertex operator algebras, generalized
doubles and dual pairs, {\itshape Math. Z.} {\bfseries 241}
(2002), 397--423.

\bibitem{FZ}
V. A. Fateev and A. B. Zamolodchikov, Conformal quantum field
theory models in two dimensions having $\Z_3$ symmetry, {\itshape
Nuclear Physics} {\bfseries B280} (1987), 644--660.

\bibitem{FHL}  I. B. Frenkel, Y. Huang and J. Lepowsky, {\itshape
On axiomatic approaches to vertex operator algebras and modules},
Mem. Amer. Math. Soc. 104, 1993.

\bibitem{FKW}
E. Frenkel, V. G. Kac, and M. Wakimoto, Characters and fusion
rules for $W$-algebras via quantized Drinfel${}^{\prime}$d-Sokolov
reduction, Comm. Math. Phys. {\bfseries 147} (1992), 295--328.

\bibitem{FLM}
I. B. Frenkel, J. Lepowsky and A. Meurman, {\itshape Vertex
Operator Algebras and the Monster}, Pure and Applied Math., Vol.
{\bfseries 134}, Academic Press, 1988.

\bibitem{K}
V. G. Kac,
Infinite-dimensional Lie algebras.
Third edition. Cambridge University Press, Cambridge, 1990.

\bibitem{KR}
V. G. Kac and A. K. Raina, {\itshape Highest Weight
representations of Infinite Dimensional Lie Algebras}, World
Scientific, 1987.

\bibitem{KW}
V. G. Kac and M. Wakimoto, Modular and conformal invariance
constraints in representation theory of affine algebras, Adv. in
Math. {\bfseries 70} (1988), 156--236.

\bibitem{KLY1}
K. Kitazume, C. Lam and H. Yamada, Decomposition of the
moonshine vertex operator algebra as Virasoro modules, {\itshape
J. Algebra}, {\bfseries 226} (2000), 893--919.

\bibitem{KLY2}
K. Kitazume, C. Lam and H. Yamada, $3$-state Potts model,
moonshine vertex operator algebra and $3A$ elements of the monster
group, to appear in {\itshape International Mathematics Research
Notices}.

\bibitem{KMY}
M. Kitazume, M. Miyamoto and H. Yamada, Ternary codes and vertex
operator algebras, {\itshape J. Algebra}, {\bfseries 223} (2000),
379--395.

\bibitem{LY}
C. Lam and H. Yamada, $\mathbb{Z}_{2}\times
\mathbb{Z}_{2}$ codes and vertex operator algebras, {\itshape J.
Algebra} {\bfseries 224} (2000), 268--291.

\bibitem{L}
J. Lepowsky, Calculus of twisted vertex operators, {\itshape Proc.
Natl. Acad. Sci. USA} {\bfseries 82} (1985), 8295--8299.

\bibitem{Li}
H. Li, Symmetric invariant bilinear forms on vertex operator
algebras, {\itshape J. Pure and Applied Alg.} {\bfseries 96}
(1994), 279--297.

\bibitem{M1}
M. Miyamoto, Griess algebras and conformal vectors
in vertex operator algebras, {\itshape J. Algebra} {\bfseries 179}
(1996), 523--548.

\bibitem{M2}
M. Miyamoto, $3$-State Potts model and automorphism of vertex
operator algebra of order $3$, {\itshape J. Algebra} {\bfseries 239}
(2001), 56--76.

\bibitem{M3}
M. Miyamoto,
A new construction of the moonshine vertex operator algebra over
the real number field, to appear in {\itshape Ann. of Math.}

\bibitem{MT}
M. Miyamoto and K. Tanabe, Uniform product of $A_{g,n}(V)$ for an
orbifold model $V$ and $G$-twisted Zhu algebra, math.QA/0112054.

\bibitem{R}
A. Rocha-Caridi, Vacuum vector representations of
the Virasoro algebra, {\itshape in} ``Vertex Operators in
Mathematics and Physics,'' Publications of the Mathematical
Sciences Research Institute, Vol. {\bfseries 3}, Springer-Verlag,
Berlin/New York, 1984, pp. 451--473.

\bibitem{W}
M. Wakimoto, Infinite-dimensional Lie algebras,
Translated from the 1999 Japanese original by Kenji Iohara,
Translations of Mathematical Monographs {\bfseries 195},
Iwanami Series in Modern Mathematics,
American Mathematical Society, Providence, RI, 2001.

\bibitem{Wa1}
W. Wang, Rationality of Virasoro vertex operator algebras,
{\itshape Duke Math. J.} {\bfseries 71}, {\itshape Inter. Math.
Res. Notice} (1993), 197--211.

\bibitem{Wa2} W. Wang,
Classification of irreducible modules of
$\W_3$ algebra with $c=-2$,
{\itshape Comm. Math. Phys.} {\bfseries 195} (1998), 113--128.

\bibitem{Z}
Y. Zhu, Modular invariance of characters of vertex operator
algebras, {\itshape J. Amer. Math. Soc.} {\bfseries 9} (1996),
237--302.

\end{thebibliography}
\end{document}